\documentclass[leqno, amstex,12pt, amssymb, final]{article}

\usepackage{mathtext}
\usepackage[cp1251]{inputenc}
\usepackage[T2A]{fontenc}
\usepackage[dvips]{graphicx}
\usepackage{amsmath}
\usepackage{amssymb}
\usepackage{amsxtra}
\usepackage{latexsym}
\usepackage{ifthen}
\date{}

\newcounter{lemma}[section]

\newcounter{corollary}[section]

\newcounter{remark}[section]

\newcounter{theorem}[section]

\newcounter{proposition}[section]

\newcounter{example}

\numberwithin{equation}{section}

\begin{document}

\markboth{M.~MATELJEVI\'{C}, E.~SEVOST'YANOV}{\centerline{On Koebe
theorem...}}

\def\cc{\setcounter{equation}{0}
\setcounter{figure}{0}\setcounter{table}{0}}

\overfullrule=0pt


\title{
{\bf On Koebe theorem for Sobolev mappings on Riemannian surfaces}}

\author{Miodrag Mateljevi\'{c}, Evgeny Sevost'yanov\footnote{Corresponding author, esevostyanov2009@gmail.com}}

\date{\today}
\maketitle

\begin{abstract}
We consider mappings between domains of two Riemannian surfaces. The
main subject of the investigation are families of homeomorphisms
with moduli inequalities. For such families, we obtain the uniform
openness (uniform lightness). As a consequence, we obtain similar
results for homeomorphisms of Sobolev classes with finite distortion
whenever its maximal dilatations have a majorant integrable on
almost all concentric circles. In particular, we obtain the
Koebe-type theorem for both direct and inverse Sobolev mappings.
\end{abstract}



\section{Introduction}

The known Koebe theorem states that for broad classes of analytic
functions, the image of the unit disk contains a disk centered at
zero with radius~1/4. More precisely, the following statement holds,
see, e.g., ~\cite[Theorem~1.3]{CG}, cf. \cite{Cr$_1$}, \cite{E},
\cite{M$_1$}--\cite{M$_3$}, \cite{Ra} and \cite{RS$_2$}.

\medskip
{\bf Theorem A.} {\it Let $f:{\Bbb D}\rightarrow {\Bbb C}$ be an
univalent analytic function such that $f(0)=0$ and
$f^{\,\prime}(0)=1.$ Then the image of ${\Bbb D}$ under the mapping
$f$ covers the open disk centered at $0$ of radius one-quarter, that
is, $f({\Bbb D})\supset B(0, 1/4).$}

\medskip
Relatively recently, we have investigated the connection between
statements of Koebe type and the distortion estimates for modulus of
families of paths, see, e.g., \cite{RS$_1$}, \cite{RS$_2$} and
\cite{ST$_1$}--\cite{ST$_2$}. The main goal of our manuscript is to
show how the tools and ideas we developed earlier work on Riemannian
surfaces. In particular, in \cite{IRS} and \cite{RS$_1$} we
established such a connection for homeomorphisms with Poletskii's
inequality of the Euclidean space whose (quasiconformality)
characteristic is integrable over concentric spheres. Below we prove
similar result for Riemannian surfaces of hyperbolic type. The
necessary tools for investigation of mappings between Riemannian
surfaces were were developed in \cite{RV$_1$}--\cite{RV$_2$}.

\medskip
The definitions concerning Riemannian surfaces may be found, for
example, in~\cite{RV$_1$} and \cite[Chapter~17]{Sev}. The
work~\cite{RV$_1$} also presents all the necessary information
concerning prime ends on Riemannian surfaces, and therefore they are
omitted below. In what follows, the Riemannian surfaces ${\Bbb S}$
and ${\Bbb S}_*$ have hyperbolic type with metrics $\widetilde{h}$
and $\widetilde{h}_*,$ respectively (see, e.g., \cite{RV$_1$},
\cite{RV$_2$}). Besides that, $ds_{\widetilde{h}}$ and
$d\widetilde{v},$ $ds_{\widetilde{h}_*}$ and $d\widetilde{v}_*$
denote the elements of length and area on the Riemannian surfaces
${\Bbb S}$ and ${\Bbb S}_*,$ respectively. We set
$$\widetilde{B}(p_0, r):=\{p\in {\Bbb S}: \widetilde{h}(p,
p_0)<r\}, \quad\widetilde{S}(p_0, r):=\{p\in {\Bbb S}:
\widetilde{h}(p, p_0)=r\}\,.$$
Below we also make extensive use of the concept of a normal
neighborhood of a point on a Riemannian surface, that is, a
neighborhood that is homeomorphic to some subset of the unit disk,
and in which the distance between points on the surface coincides
with the hyperbolic distance between the corresponding points in the
unit disk, see, e.g., \cite[Proposition~18.1 and Lemma~18.3]{Sev}.
In more detail, a neighborhood $U$ of a point $x_0\in {\Bbb S}$ is
{\it normal,} if there exists a hyperbolic isometry $\varphi$ of
this neighborhood onto some neighborhood $V$ lying in the unit disk
${\Bbb D}\subset{\Bbb C},$ while $\varphi(x_0)=0,$ $\widetilde{h}(x,
y)=h(\varphi(x), \varphi(y))$ for any $x, y\in U,$ and $h(x, y)$ is
a hyperbolic metric,
\begin{equation}\label{eq7AC}h(z_1, z_2)=\log\,\frac{1+t}{1-t}\,,\quad
t=\frac{|z_1-z_2|}{|1-z_1\overline{z_2}|}\,.
\end{equation}
For the existing of the normal neighborhoods see, e.g.,
\cite[Proposition~18.1, Lemma~18.3]{Sev}.

\medskip
As usual, a path $\gamma$ on the Riemannian surface ${\Bbb S}$ is
defined as a continuous mapping $\gamma:I\rightarrow {\Bbb S},$
where $I$ is a finite segment, an interval or a half-interval of a
real axis. Let $\Gamma$ be a family of paths $\gamma$ in ${\Bbb S}.$
A Borel function $\rho:{\Bbb S}\rightarrow [0, \infty]$ is called
{\it admissible for $\Gamma,$} if
$\int\limits_{\gamma}\rho(p)\,ds_{\widetilde{h}}(p)\geqslant 1$ for
any (locally rectifiable) path $\gamma \in \Gamma.$ The latter is
briefly written in the form: $\rho\in {\rm adm}\,\Gamma.$  A {\it
modulus} of the family $\Gamma$ is defined as
$$M(\Gamma):=\inf\limits_{\rho\in {\rm adm}\,\Gamma}\int\limits_{\Bbb
S}\rho^2(p)\,d\widetilde{v}(p)\,.$$
To avoid confusion, we introduce separately the notation for the
modulus of the family $\Gamma$ in the Euclidean sense, putting
\begin{equation}\label{eq1C}
M_e(\Gamma):=\inf\limits_{\rho\in {\rm
adm}_e\,\Gamma}\int\limits_{\Bbb D}\rho^2(z)\,dm(z)\,,
\end{equation}
where $\rho\in {\rm adm}_e\,\Gamma$ if and only if
$\int\limits_{\gamma}\rho(z)\,|dz|\geqslant 1$ for any (locally
rectifiable) path $\gamma \in \Gamma.$ Finally,
\begin{equation}\label{eq1C_1}
M_h(\Gamma):=\inf\limits_{\rho\in {\rm
adm}_h\,\Gamma}\int\limits_{\Bbb D}\rho^2(z)\,dv(z)\,,
\end{equation}
where $dv(z)$ is an element of the hyperbolic measure,
$$dv(z)=\frac{4\,dm(z)}{(1-|z|^2)^2}\,,$$
and $\rho\in {\rm adm}_h\,\Gamma$ if and only if
$\int\limits_{\gamma}\rho(z)\,ds_{h}(z)\geqslant 1$ for any (locally
rectifiable) path $\gamma \in \Gamma,$ where $h$ is a hyperbolic
metric, see~(\ref{eq7AC}). Observe that
\begin{equation}\label{eq10}
M_h(\Gamma)=M_e(\Gamma)
\end{equation}
whenever $\Gamma$ is a family of paths in the unit disk, see
\cite[relation~(9.31)]{Sev}.
Let $Q:{\Bbb S}\rightarrow [0, \infty]$ be a function measurable
with respect to the measure $\widetilde{v},$ $Q(x)\equiv 0$ for
$x\not\in D\subset {\Bbb S}.$ We say that $f:D\rightarrow {\Bbb
S}_*$ is a {\it ring $Q$-mapping at $p_0\in D,$ }if the relation
\begin{equation}\label{eq3*!!A}
M(f(\Gamma(\widetilde{S}(p_0, r_1),\,\widetilde{S}(p_0,
r_2),\,D)))\leqslant \int\limits_A Q(p)\cdot \eta^2(\widetilde{h}(p,
p_0))\,d\widetilde{v}(p)
\end{equation}
holds for some $r_0=r(p_0)>0,$ any ring
$\widetilde{A}=\widetilde{A}(p_0, r_1, r_2)=\{p\in {\Bbb S}:
r_1<\widetilde{h}(p, p_0)<r_2\},$ any $0<r_1<r_2<r_0$ and all
nonnegative Lebesgue measurable function $\eta:(r_1,r_2)\rightarrow
[0,\infty ]$ that
\begin{equation}\label{eq8BC}
\int\limits_{r_1}^{r_2}\eta(r)\,dr\geqslant 1\,.
\end{equation}
holds.

\medskip
Let $D$ be a domain in ${\Bbb S}$ and let $K$ be a compact subset of
$D.$ A family $\frak{F}$ of mappings $f:D\rightarrow {\Bbb S}_*$ is
called {\it uniformly light on $K,$} if the following condition
holds: given $\varepsilon>0$ there is $\delta_1(\varepsilon, K)>0$
such that $\widetilde{h}_*(f(C))\geqslant \delta_1$ for any $f\in
\frak{F}$ and any continuum $C\subset K$ with
$\widetilde{h}(C)\geqslant \varepsilon.$ A family $\frak{F}$ of
mappings $f:D\rightarrow {\Bbb S}_*$ is called {\it equi-uniform on
$K,$} if the following condition holds: given $\varepsilon>0$ there
is $\delta_1(\varepsilon, K)>0$ such that $\widetilde{h}_*(f(x),
f(y))\geqslant \delta_1$ for any $f\in \frak{F}$ and any $x, y\in K$
with $\widetilde{h}(x, y)\geqslant \varepsilon.$ Obviously, if a
family is equi-uniform on $K,$ then it also is uniformly light
on~$K.$

\medskip
A family $\frak{F}$ of mappings $f:D\rightarrow {\Bbb S}_*$ is
called {\it uniformly open at a point $x_0\in D,$} if for every
$\varepsilon_0>0,$ $\varepsilon_0<\widetilde{h}(x_0, \partial D),$
there is $r_0=r_0(x_0, \varepsilon_0)>0$ such that
$\widetilde{B}(f(x_0), r_0)\subset f(\widetilde{B}(x_0,
\varepsilon_0))$ for every $f\in \frak{F}.$ A family $\frak{F}$ of
mappings $f:D\rightarrow {\Bbb S}_*$ is called {\it uniformly open
on $K,$} if for every $\varepsilon_0>0,$
$\varepsilon_0<\widetilde{h}(K, \partial D),$ there exists
$r_0=r_0(K, \varepsilon_0)>0$ such that $\widetilde{B}(f(x_0),
r_0)\subset f(\widetilde{B}(x_0, \varepsilon_0))$ for every $f\in
\frak{F}$ and every $\widetilde{B}(x_0, \varepsilon_0)\subset K.$

\medskip
Let $(X, d)$ be a metric space, $A, B \subset X,$
$d(A)=\sup\limits_{x, y\in A}d(x, y),$ $d(A, B)=\inf\limits_{x\in A,
y\in B}d(x, y).$ For each real number $n\geqslant 1,$ we define {\it
the Loewner function} $\Phi_n:(0, \infty)\rightarrow [0, \infty)$ on
$X$ as
\begin{equation}\label{eq2H}
\Phi_n(t)=\inf\{M_n(\Gamma(E, F, X)): \Delta(E, F)\leqslant t\}\,,
\end{equation}
where $M_n$ is the $n$-modulus of family of paths in $X$ (see
\cite{He}), besides that, the infimum is taken over all disjoint
nondegenerate continua $E$ and $F$ in $X$ and
$$\Delta(E, F):=\frac{d(E,
F)}{\min\{d(E), d(F)\}}\,.$$
A pathwise connected metric measure space $(X, \mu)$ is said to be a
{\it Loewner space} of exponent $n,$ or an $n$-Loewner space, if the
Loewner function $\Phi_n(t)$ is positive for all $t> 0$ (see
\cite[Section~2.5]{MRSY} or \cite[Ch.~8]{He}). Observe that, ${\Bbb
R}^n$ is a Loewner space (see \cite[Theorem~8.2]{He}).

\medskip
Let $X$ be a Loewner space in which the relation $\mu(B_R)\leqslant
C^{\,*}R^n$ holds for some constant $C^{\,*}\geqslant 1,$ for some
exponent $n>0$ and for all closed balls $B_R$ of radius $R>0.$ Let
$\Phi_n(t)$ be a Loewner function in~(\ref{eq2H}), corresponding to
the space $X.$ Then, by virtue of~\cite[Theorem~8.23]{He}, there is
$\widetilde{\delta}_0>0$ and some constant $C>0$ such that
\begin{equation}\label{eq3C}
\Phi_n(t)\geqslant C\log\frac{1}{t}
\end{equation}
for $0<t<\widetilde{\delta}_0.$ It is well-known that, ${\Bbb R}^n$
is a Loewner space, see, e.g., \cite[Theorem~8.2]{He}. In
particular, ${\Bbb C}\sim {\Bbb R}^2$ is a Loewner space, as well.
Obviously, the above condition $\mu(B_R)\leqslant C^{\,*}R^n$ holds
in ${\Bbb C}$ with the Lebesgue measure $\mu=m,$ $n=2$ and
$C^{\,*}=\pi.$ Thus, the relation~(\ref{eq3C}) obviously holds in
${\Bbb C}$ with some constant $C>0.$

\medskip
Given domains $D$ and $D_*$ in ${\Bbb S}$ and ${\Bbb S}_*,$
respectively, a function $Q:D\rightarrow [0, \infty]$ measurable by
the measure $\widetilde{v},$ a compact set $K\subset D$ and a number
$\delta>0$ denote by $\frak{F}^{\delta}_{K, Q}(D, D_*)$ a family of
homeomorphisms $f:D\rightarrow D_*$ such that, for any $p_0\in D$
there is $r_0=r_0(p_0)>0,$ $r_0<\sup\limits_{p\in D}\widetilde{h}(p,
p_0),$ such that the relations~(\ref{eq3*!!A})--(\ref{eq8BC}) hold
for any $0<r_1<r_2<r_0$ and, besides that, $\widetilde{h}_*(f(K),
\partial f(D))=\inf\limits_{x\in f(K), y\in \partial f(D)}
\widetilde{h}_*(x, y)\geqslant \delta.$ The following results hold.

\begin{theorem}\label{th1}
{\,Let $D$ and $D_*$ be domains in ${\Bbb S}$ and ${\Bbb S}_*,$
respectively, while $D$ has at least two different boundary points
$a$ and $b$ and, besides that, $\overline{D_*}$ is a compactum in
${\Bbb S}_*.$ Assume that, for each point $p_0\in D$ there is
$r_0=r_0(p_0)>0,$ $r_0<\sup\limits_{p\in D}\widetilde{h}(p, p_0),$
such that for any $0<r_1<r_2<r_0$ there is a set $E_1\subset[r_1,
r_2]$ of a positive linear Lebesgue measure such that the function
$Q$ is integrable over $\widetilde{S}(p_0, r)$ for every $r\in E_1.$
Then the family $\frak{F}^{\delta}_{K, Q}(D, D_*)$ is uniformly
light on $K;$ moreover, $\frak{F}^{\delta}_{K, Q}(D, D_*)$ is
equi-uniform on $K.$

\medskip
Moreover, if (\ref{eq3*!!A})--(\ref{eq8BC}) hold for any $p_0\in D,$
any $0<r_1<r_2<\sup\limits_{p\in D}\widetilde{h}(p, p_0),$ any
$f\in\frak{F}^{\delta}_{K, Q}(D, D_*)$ and $Q\in L^1(D),$ then there
exist $C_1>0$ such that
\begin{equation}\label{eq1D}
\widetilde{h}_*(f(x), f(y))\geqslant
C_1\cdot\exp\left\{-\frac{2\Vert Q\Vert_1}{C{(\widetilde{h}(x,
y))}^2}\right\}
\end{equation}
for all $x, y\in K,$ $x\ne y,$ and every $f\in \frak{F}^{\delta}_{K,
Q}(D, D_*),$ where $C$ is a constant from (\ref{eq3C}) for $n=2$ and
$X={\Bbb C}.$ }
\end{theorem}

\begin{corollary}\label{cor1}
{\, The statement of Theorem~\ref{th1} remains true if we replace
the assumption on $Q$ mentioned above concerning the integrability
of $Q$ over circles $\widetilde{S}(p_0, r)$ by a simpler condition:
$Q\in L_{\rm loc}^1(D).$}
\end{corollary}

\begin{corollary}\label{cor2}{\bf \,(An analogue of Koebe-Bloch
theorem)}.~{Under assumptions of Theorem~\ref{th1}, the family
$\frak{F}^{\delta}_{K, Q}(D, D_*)$ is uniformly open on $K.$ In
particular, the latter is true if we replace the assumption on $Q$
mentioned in Theorem~\ref{th1} concerning the integrability of $Q$
over circles $\widetilde{S}(p_0, r)$ by a simpler condition: $Q\in
L_{\rm loc}^1(D).$}
\end{corollary}

Note that, the condition $\widetilde{h}(f(K),
\partial f(D))\geqslant \delta$ in
$\frak{F}^{\delta}_{K, Q}(D, D_*)$ is somewhat difficult to verify,
so we consider now some different class of mappings used another
condition.

\medskip
Given domains $D, D_*$ in ${\Bbb S}$ and ${\Bbb S}_*,$ respectively,
a non-degenerate continuum $A\subset D,$ a measurable function
$Q:D\rightarrow [0, \infty]$ and a number $\delta>0$ we denote by
$\frak{B}^{\delta}_{A, Q}(D, D_*)$ a family of all homeomorphisms
$f:D\rightarrow D_*$ satisfying~(\ref{eq3*!!A})--(\ref{eq8BC}) for
all $p_0\in D,$ some $r_0=r_0(p_0)>0,$ $r_0<\sup\limits_{p\in
D}\widetilde{h}(p, p_0),$ and all $0<r_1<r_2<r_0$ such that
$\widetilde{h}_*(f(A))\geqslant \delta.$

\begin{theorem}\label{th2}
{\,Let $D$ and $D_*$ be domains in ${\Bbb S}$ and ${\Bbb S}_*,$
respectively, while $D$ has at least two different boundary points
$a$ and $b$ and, besides that, $\overline{D_*}$ is a compactum in
${\Bbb S}_*.$ Assume that,

for each point $p_0\in D$ there is $r_0=r_0(p_0)>0,$
$r_0<\sup\limits_{p\in D}\widetilde{h}(p, p_0),$ such that for any
$0<r_1<r_2<r_0$ there is a set $E_1\subset[r_1, r_2]$ of a positive
linear Lebesgue measure such that the function $Q$ is integrable
over $\widetilde{S}(p_0, r)$ for every $r\in E_1.$ Then the family
$\frak{B}^{\delta}_{A, Q}(D, D_*)$ is uniformly light on every
compactum $K\subset D$ moreover, $\frak{B}^{\delta}_{A, Q}(D, D_*)$
is equi-uniform on $K.$

\medskip
Moreover, if (\ref{eq3*!!A})--(\ref{eq8BC}) hold for any $p_0\in D,$
any $0<r_1<r_2<\sup\limits_{p\in D}\widetilde{h}(p, p_0),$ any $f\in
\frak{B}^{\delta}_{A, Q}(D, D_*)$ and $Q\in L^1(D),$ then there
exist $C_1>0$ such that the relation~(\ref{eq1D}) holds
for all $x, y\in K,$ $x\ne y,$ and every $f\in \frak{B}^{\delta}_{A,
Q}(D, D_*),$ where $C$ is a constant from (\ref{eq3C}) for $n=2$ and
$X={\Bbb C}.$ }
\end{theorem}

\begin{corollary}\label{cor6}
{\, The statement of Theorem~\ref{th2} remains true if we replace
the assumption on $Q$ concerning the integrability of $Q$ over
circles $\widetilde{S}(p_0, r)$ by a simpler condition: $Q\in L_{\rm
loc}^1(D).$}
\end{corollary}

\begin{corollary}\label{cor7}{\bf \,(An analogue of Koebe-Bloch
theorem)}.~{ Under assumptions of Theorem~\ref{th2}, the family
$\frak{B}^{\delta}_{A, Q}(D, D_*)$ is uniformly open on every
compactum $K\subset D.$ In particular, the latter is true if we
replace the assumption on $Q$ mentioned in Theorem~\ref{th2}
concerning the integrability of $Q$ over circles $\widetilde{S}(p_0,
r)$ by a simpler condition: $Q\in L_{\rm loc}^1(D).$}
\end{corollary}

\section{On the approaching continuum}

The following statement holds.

\begin{lemma}\label{lem3} {\bf(V\"{a}is\"{a}l\"{a}'s lemma).}
{\, Let $D$ be a domain in $\overline{\Bbb C},$ and let $x_0\in D.$
Then for each $P>0$ and each neighborhood $U$ of point $x_0$ there
is a neighborhood $V\subset U$ of the same point such that
$M_e(\Gamma(E, F, D))> P$ for any continua $E, F \subset D $
intersecting $\partial U$ and $\partial V.$}
\end{lemma}

The proof of Lemma~\ref{lem3} is essentially given by
V\"{a}is\"{a}l\"{a} in~\cite[(10.11)]{Va}, however, we also given a
formal proof (\cite[Lemma~2.2]{SevSkv$_1$}).~$\Box$

\medskip
Given domains $D$ and $D_*$ in ${\Bbb S}$ and ${\Bbb S}_*,$
respectively, a non-degenerate continuum $A\subset D,$ a function
$Q:D\rightarrow [0, \infty]$ measurable with respect to the measure
$v$ on ${\Bbb S}$ and numbers $\delta, r>0$ we denote by
$\frak{A}^{\delta, r}_{A, Q}(D, D_*)$ a family of all homeomorphisms
$f:D\rightarrow D_*$ such that, for any point $p_0\in D$ there is
$r_0=r_0(p_0)>0,$ $r_0<\sup\limits_{p\in D}\widetilde{h}(p, p_0),$
such that the relations~(\ref{eq3*!!A})--(\ref{eq8BC}) hold for any
$0<r_1<r_2<r_0,$ besides that, $\widetilde{h}_*(f(A))\geqslant
\delta$ and the following condition holds:
$\widetilde{h}_*(E)\geqslant r$ whenever $E$ is a component of
$\partial f(D).$

\medskip
We have proven different versions of the following lemma, see, e.g.,
\cite[Lemma~3.2]{RS$_1$}, \cite[Lemma~4.1]{SevSkv$_2$}.

\begin{lemma}\label{lem2}{\, Let $D$ and $D_*$ be domains in ${\Bbb S}$ and ${\Bbb S}_*,$
respectively, while $\overline{D_*}$ is a compactum in ${\Bbb S}_*.$
Let $A$ be a non-degenerate continuum in $D$ and let
$f_m:D\rightarrow D_*$ be a sequence of homeomorphisms in
$\frak{A}^{\delta, r}_{A, Q}(D)$ for some function~$Q,$ some
continuum $A\subset D$ and numbers $\delta, r>0.$ Assume that, for
each point $p_0\in D$ there is $r_0=r_0(p_0)>0,$
$r_0<\sup\limits_{p\in D}\widetilde{h}(p, p_0),$ such that for any
$0<r_1<r_2<r_0$ there is a set $E_1\subset[r_1, r_2]$ of a positive
linear Lebesgue measure such that the function $Q$ is integrable
with respect to $s_h$ over the spheres $\widetilde{S}(x_0, r)$ for
every $r\in E_1.$ If $\partial f_m(D)\ne \varnothing,$ then there
exists $\delta_1>0$ such that $\widetilde{h}_*(f_m(A),
\partial f_m(D))\geqslant \delta_1$ for any $m=1,2,\ldots .$ }
\end{lemma}

\begin{proof}
We prove Lemma~\ref{lem2} by the contradiction, cf. Lemma~4.1 in
\cite{SevSkv$_2$}. Assuming that the above conclusion is not true,
for each $k\in {\Bbb N}$ we may find some $m=m_k\in {\Bbb N}$ such
that $\widetilde{h}_*(f_{m_k}(A),
\partial f_{m_k}(D))<1/k.$ Due to renumbering, if required, we may
write $f_m$ instead of $f_{m_k},$ i.e., $\widetilde{h}_*(f_{m}(A),
\partial f_m(D))<1/m,$ $m=1,2,\ldots. $ Since $\overline{D_*}$ is a compactum in
${\Bbb S}_*$, $f_{m}(A)$ is a compactum as a continuous image of a
compactum $A\subset D$ under $f_{m},$ we may find $x_m\in f_{m}(A)$
and $y_m\in \partial f_m(D)$ such that $\widetilde{h}_*(f_{m}(A),
\partial f_m(D))=\widetilde{h}_*(x_m, y_m)<1/m.$ Let $E_m$ be a component of $\partial
f_m(D)$ containing $y_m.$ Due to the compactness of
$\overline{D_*},$ we may assume that $y_m\rightarrow y_0$ as
$m\rightarrow \infty;$ then also $x_m\rightarrow y_0$ as
$m\rightarrow \infty.$

\medskip
Put $P>0$ and $U=\widetilde{B}(y_0, r_0)=\{y\in {\Bbb S}_*:
\widetilde{h}_*(y, y_0)<r_0\},$ where $2r_0:=\min\{r/2, \delta/2\},$
and $r$ and $\delta$ are numbers from the definition of the class
$\frak{A}^{\delta, r}_{A, Q}(D).$ We may consider that $U$ lies in
the normal neighborhood $U_0$ of $y_0.$ Now, there is a mapping
$\varphi:U_0\rightarrow {\Bbb C}$ such that $\varphi$ maps
isometrically $\widetilde{B}(y_0, r_0)$ onto $B_h(0, r_0)=\{y\in
{\Bbb D}: h(y, 0)<r_0\}.$ Let $E^{\,*}_m$ be a component of $E_m\cap
U$ containing $y_m.$ Observe that, $E^{\,*}_m$ is closed (see
\cite[Theorem~2.II.4.41]{Ku$_2$}), besides that, $E_m\ne E^{\,*}_m$
for sufficiently large $m\in{\Bbb N},$ because $y_m\rightarrow y_0$
as $m\rightarrow \infty,$ $y_m\in E_m,$
$\widetilde{h}_*(E^{\,*}_m)\leqslant
\widetilde{h}_*(\widetilde{B}(y_0, r_0))\leqslant 2r_0\leqslant
\delta/2$ and $\widetilde{h}_*(E_m)\geqslant \delta.$ Now, by
Janiszewski's theorem, $E^{\,*}_m\cap \partial U\ne\varnothing.$ Let
$F_m$ be a component of $f_m(A)\cap U$ containing $x_m.$  Arguing
similarly, we may show that $F_m$ is closed, $F_m\ne f_m(A)$ for
sufficiently large $m\in{\Bbb N}$ and $F_m\cap \partial
U\ne\varnothing.$ Let us consider the neighborhood $V_*$ of the
origin corresponding to given $P>0$ and the neighborhood
$\varphi(U)$ in ${\Bbb C}$ in the context of Lemma~\ref{lem3}. We
put $V:=\varphi^{\,-1}(V_*).$ Arguing similarly to above, we may
show that $E^{\,*}_m\cap
\partial V\ne\varnothing\ne F_m\cap
\partial V.$ Now, $\varphi(E^{\,*}_m)\cap
\partial V_*\ne\varnothing\ne \varphi(F_m)\cap
\partial V_*$ and $\varphi(E^{\,*}_m)\cap
\partial \varphi(U)\ne\varnothing\ne \varphi(F_m)\cap
\partial \varphi(U).$ By Lemma~\ref{lem3},
\begin{equation}\label{eq11}
M_e(\Gamma(\varphi(F_m), \varphi(E^{\,*}_m), \varphi(U)))>P\,,
\end{equation}
where $M_e$ is defined in~(\ref{eq1C}). Since the Euclidean and the
hyperbolic moduli coincide, see, e.g., \cite[relation~(9.31)]{Sev},
we have from here that
\begin{equation}\label{eq9B}
M_h(\Gamma(\varphi(F_m), \varphi(E^{\,*}_m),
\varphi(U)))=M(\Gamma(F_m, E^{\,*}_m, U))>P
\end{equation}
for sufficiently large $m=1,2,\ldots .$ We now prove that the
relation~(\ref{eq9B}) contradicts the definition of $f_m$
in~(\ref{eq3*!!A})--(\ref{eq8BC}). Indeed, let $\gamma:[0,
1]\rightarrow {\Bbb S}_*$ be a path in $\Gamma(F_m, E^{\,*}_m, {\Bbb
S}_*),$ i.e., $\gamma(0)\in F_m,$ $\gamma(1)\in E^{\,*}_m$ and
$\gamma(t)\in {\Bbb S}_*$ for $t\in (0, 1).$ Let
$t_m=\sup\limits_{\gamma(t)\in f_m(D)}t$ and let
$\alpha_m(t)=\gamma|_{[0, t_m)}(t).$ Let $\Gamma_m$ consists of all
such paths $\alpha_m,$ now $\Gamma(f_m(A), E_m, {\Bbb
S}_*)>\Gamma_m$ and by the minorization principle of the modulus
(see~\cite[Theorem~1]{Fu})
\begin{equation}\label{eq12}
M(\Gamma_m)\geqslant M(\Gamma(F_m, E^{\,*}_m, {\Bbb S}_*))\,.
\end{equation}
Let $\Delta_m:=f^{\,-1}_m(\Gamma_m).$ Now
\begin{equation}\label{eq13}
M(f_m(\Delta_m))=M(\Gamma_m)\,.
\end{equation}
Observe that, $\beta_m\in \Delta_m$ if and only if
$\beta_m(t)=f^{\,-1}_m(\alpha_m(t))=f^{\,-1}_m(\gamma|_{[0,
t_m)}(t))$ for some $\gamma\in \Gamma(F_m, E^{\,*}_m, {\Bbb S}_*),$
$t_m=\sup\limits_{\gamma(t)\in f_m(D)}t.$ Since homeomorphisms
preserve the boundary, $\beta_m(t)\rightarrow \partial D$ as
$t\rightarrow t_m.$ In particular, it follows that $\partial
D\ne\varnothing.$

Let $\varepsilon:=\widetilde{h}\,(A, \partial D).$ Given $x\in {\Bbb
S},$ let $\varepsilon_x$ be the number such that
$0<\varepsilon_x<\min\{\varepsilon, r_0(x)\},$ and, besides that,
the disk $\overline{\widetilde{B}(x, \varepsilon_x)}$ lie in some a
normal neighborhood $U_x$ of $x.$ We cover the continuum $A$ with
balls $\widetilde{B}(x, \varepsilon_x),$ $x\in A.$ Since $A$ is a
compactum, we may assume that $A\subset
\bigcup\limits_{i=1}^{M_0}\widetilde{B}(x_i, \varepsilon_i/4),$
$x_i\in A,$ $\varepsilon_i=\varepsilon_{x_i},$ $i=1,2,\ldots, M_0,$
$1\leqslant M_0<\infty.$ By the definition, $M_0$ depends only on
$A;$ in particular, $M_0$ does not depend on $m.$
Observe that,
\begin{equation}\label{eq6D}
\Delta_m=\bigcup\limits_{i=1}^{M_0}\Gamma_{mi}\,,
\end{equation}
where $\Gamma_{mi}$ consists of all paths $\gamma:[0,
t_m)\rightarrow D$ in $\Delta_m$ such that $\gamma(0)\in
\widetilde{B}(x_i, \varepsilon/4)$ and $\gamma(t)\rightarrow
\partial D$ as $t\rightarrow t_m-0.$ We now show that
\begin{equation}\label{eq7C}
\Gamma_{mi}>\Gamma(\widetilde{S}(x_i, \varepsilon_i/4),
\widetilde{S}(x_i, \varepsilon_i/2), \widetilde{A}(x_i,
\varepsilon_i/4, \varepsilon_i/2))\,.
\end{equation}
Indeed, let $\gamma\in \Gamma_{mi},$ in other words, $\gamma:[0,
t_m)\rightarrow D,$ $\gamma\in \Delta_m,$ $\gamma(0)\in
\widetilde{B}(x_i, \varepsilon_i/4)$ and $\gamma(t)\rightarrow
\partial D$ as $t\rightarrow t_m-0.$ Now, by the definition of
$\varepsilon,$ $|\gamma|\cap \widetilde{B}(x_i,
\varepsilon_i/4)\ne\varnothing\ne |\gamma|\cap (D\setminus
\widetilde{B}(x_i, \varepsilon_i/4)).$ Therefore, by
\cite[Theorem~1.I.5.46]{Ku$_2$} there is $0<t^{\,*}_1<t_m$ such that
$\gamma(t^{\,*}_1)\in \widetilde{S}(x_i, \varepsilon_i/4).$ We may
assume that $\gamma(t)\not\in \widetilde{B}(x_i, \varepsilon_i/4)$
for $t>t^{\,*}_1.$ Put $\gamma_1:=\gamma|_{[t^{\,*}_1, t_m]}.$
Similarly, $|\gamma_1|\cap \widetilde{B}(x_i,
\varepsilon_i/2)\ne\varnothing\ne |\gamma_1|\cap (D\setminus
\widetilde{B}(x_i, \varepsilon_i/2)).$ By
\cite[Theorem~1.I.5.46]{Ku$_2$} there is $t^{\,*}_1<t^{\,*}_2<t_m$
with $\gamma(t^{\,*}_2)\in \widetilde{S}(x_i, \varepsilon_i/2).$ We
may assume that $\gamma(t)\in \widetilde{B}(x_i, \varepsilon_i/2)$
for $t<t^{\,*}_2.$ Put $\gamma_2:=\gamma|_{[t^{\,*}_1, t^{\,*}_2]}.$
Then, the path $\gamma_2$ is a subpath of $\gamma,$ which belongs to
the family $\Gamma(\widetilde{S}(x_i, \varepsilon_i/4),
\widetilde{S}(x_i, \varepsilon_i/2), \widetilde{A}(x_i,
\varepsilon_i/4, \varepsilon_i/2)).$ Thus, the relation~(\ref{eq7C})
is established. By Fubini's theorem for Riemannian surfaces (see,
e.g., \cite[Lemma~17.1]{Sev}),
\begin{equation}\label{eq8C}
M(f_{m}(\Gamma(\widetilde{S}(x_i, \varepsilon_i/4),
\widetilde{S}(x_i, \varepsilon_i/2)), \widetilde{A}(x_i,
\varepsilon_i/4, \varepsilon_i/2)))\leqslant \frac{1}{I_i}<\infty\,,
\end{equation}
where $I_i=I_i(x_i, \varepsilon_i/4,
\varepsilon_i/2)=\int\limits_{\varepsilon_i/4}^{\varepsilon_i/2}
\frac{dr}{rq_{x_i}(r)}$ and $I_i<\infty$ because $q_{x_i}(r)<\infty$
for $r\in E_1$ and some set of positive linear Lebesgue measure
$E_1.$ Finally, by~(\ref{eq12}), (\ref{eq13}), (\ref{eq6D}),
(\ref{eq7C}) and (\ref{eq8C}) we obtain that
\begin{gather}M(\Gamma(F_m, E^{\,*}_m, U))\leqslant M(\Gamma(F_m, E^{\,*}_m, {\Bbb S}_*))
\leqslant M(\Gamma_m)\nonumber\\
= M(f_m(\Delta_m))\leqslant
\sum\limits_{i=1}^{M_0}M(f_m(\Gamma_{mi}))\label{eq16}\\
\nonumber
\leqslant\sum\limits_{i=1}^{M_0}M(f_m(\Gamma(\widetilde{S}(x_i,
\varepsilon_i/4), \widetilde{S}(x_i, \varepsilon_i/2),
\widetilde{A}(x_i, \varepsilon_i/4, \varepsilon_i/2))))\leqslant
\sum\limits_{i=1}^{M_0}\frac{1}{I_i}:=C<\infty\,.
\end{gather}
Since $P$ in~(\ref{eq9B}) may be done arbitrary large, the
relations~(\ref{eq9B}) and~(\ref{eq16}) contradict each other. The
proof is complete.
\end{proof}

\begin{corollary}\label{cor4}
{\, The statement of Lemma~\ref{lem2} remains true if, instead of
the condition regarding the integrability of the function $Q$ over
spheres with respect to some set $E_1$ is replaced by a simpler
condition: $Q\in L_{\rm loc}^1(D).$}
\end{corollary}

\begin{proof} Let $r_0=r_0(p_0)>0$ be the number from the conditions of
Lemma~\ref{lem2}, defined for each $p_0\in D.$ We may consider that
$\widetilde{B}(p_0, r_0)$ lies in some normal neighborhood $U$ of
$p_0.$ Let $0<\varepsilon_1<\varepsilon_2<r_0.$ Observe that, the
Fubini theorem holds locally on Riemannian surfaces, namely,
\begin{equation}\label{eq7A}
\int\limits_{\varepsilon_1<\widetilde{h}(x,
p_0)<\varepsilon_2}Q(x)\,d\widetilde{v}(x)=
\int\limits_{\varepsilon_1}^{\varepsilon_2}\int\limits_{\widetilde{S}(p_0,
r)}Q(x)\,ds_{\widetilde{h}}(x)\,dr
\end{equation}
(see, e.g., \cite[Lemma~17.1]{Sev}). Since $Q\in L_{\rm loc}^1(D),$
(\ref{eq7A}) implies that $Q$ is integrable over $\widetilde{S}(p_0,
r)$ for almost every $r\in (\varepsilon_1, \varepsilon_2).$ In this
case, the desired conclusion follows directly from Lemma~\ref{lem2}.
\end{proof}

\begin{proposition}\label{pr4}
{\, Assume that all the conditions of Lemma~\ref{lem2} hold. Then
for any continuum $A^{\,*}$ in $D$ there is $\delta^{\,*}_1>0$ such
that $\widetilde{h}_*(f_m(A^{\,*}),
\partial f_m(D))\geqslant \delta^{\,*}_1$ for any $m=1,2,\ldots .$ In other words,
$\frak{A}^{\delta, r}_{A, Q}(D, D_*)\subset
\frak{F}^{\delta^{\,*}_1}_{A^{\,*}, Q}(D, D_*).$}
\end{proposition}

\begin{proof}
Indeed, by Lemma~\ref{lem2} there exists $\delta_1>0$ such that
$\widetilde{h}_*(f_m(A),
\partial f_m(D))\geqslant \delta_1$ for any $m=1,2,\ldots ,$ where
$A$ is a continuum from the definition of the class
$\frak{A}^{\delta, r}_{A, Q}(D).$ Assume the contrary, then we have
that $\widetilde{h}_*(f_{m_k}(A^{\,*}),
\partial f_{m_k}(D))\rightarrow 0$ for some increasing sequence of
numbers $m_k,$ $k=1,2,\ldots.$ We now join $A$ and $A^{\,*}$ by a
path $C$ in $D$ and consider the continuum $E:=A\cup A^{\,*}\cup C.$
We have $\widetilde{h}_*(f_{m_k}(E))\geqslant
\widetilde{h}_*(f_{m_k}(A))\geqslant \delta$ for $k=1,2,\ldots .$
Thus, by Lemma~\ref{lem2} there exists $\delta_2>0$ with
$\widetilde{h}_*(f_{m_k}(E),
\partial f_{m_k}(D))\geqslant \delta_2$ for any $k=1,2,\ldots .$
However, $\widetilde{h}_*(f_{m_k}(E),
\partial f_{m_k}(D))\leqslant \widetilde{h}_*(f_{m_k}(A^{\,*}),
\partial f_{m_k}(D))\rightarrow 0$ as $k\rightarrow\infty,$ because
$A^{\,*}\subset E.$ The obtained contradiction disproves the
assumption made above.
\end{proof}

The following statement holds, cf. \cite[Proposition~3.1]{RS$_1$}.

\begin{proposition}\label{pr6}
{\, Let $G_*$ and $D_*$ be domains in ${\Bbb S}_*$ such that
$G_*\subset D_*,$ and let $K_m,$ $m=1,2,\ldots,$ be a sequence of
continua in $G_*$ such that $\widetilde{h}_*(K_m, \partial
G_*)\geqslant \delta_1$ for any $m=1,2,\ldots$ and some
$\delta_1>0.$ If $\overline{D_*}$ is a compactum in ${\Bbb S}_*,$
then there is $\delta_2>0$ such that $\widetilde{h}_*(K_m, \partial
D_*)\geqslant \delta_2$ for any $m=1,2,\ldots .$ }
\end{proposition}

\begin{proof}
Assume the contrary. Now, for any $k=1,2,\ldots $ there is $m_k\in
{\Bbb N}$ such that $\widetilde{h}_*(K_{m_k}, \partial D_*)<1/k.$
Going to renumbering, if required, we may assume that
$\widetilde{h}_*(K_{m}, \partial D_*)<1/m,$ $m=1,2,\ldots.$ Since
$K_{m}$ and $\partial D_*$ are compacta, $\widetilde{h}_*(K_{m},
\partial D_*)=\widetilde{h}_*(x_m, y_m)$ for some $x_m\in K_{m}$ and
$y_m\in \partial D_*.$ Since $\overline{D_*}$ is a compactum in
${\Bbb S}_*$ and $x_m, y_m\in \overline{D_*},$ we may assume that
$x_m, y_m\rightarrow x_0$ as $m\rightarrow\infty.$ Since $\partial
D_*$ is closed, $x_0\in \partial D_*.$ On the other hand, since
$x_m\in K_m\subset G_*,$ we have that $x_0\in \overline{G_*}.$ We
now show that, the case $x_0\in \partial G_*$ is impossible. Indeed,
let $\widetilde{h}_*(x_m,
\partial G_*)=\widetilde{h}_*(x_m, z_m)$ for some $z_m\in \partial
G_*.$ Now,
$$\widetilde{h}_*(x_m,
\partial G_*)=\widetilde{h}_*(x_m, z_m)\geqslant \widetilde{h}_*(K_m,
\partial G_*)\geqslant \delta_1\,.$$
In this case, by the triangle inequality, for $z\in \partial G_*$
and sufficiently large $m=1,2,\ldots ,$
$$\widetilde{h}_*(x_0,
z)\geqslant \widetilde{h}_*(z, x_m)-\widetilde{h}_*(x_m,
x_0)\geqslant \delta_1-\frac{1}{2}\delta_1=\frac{1}{2}\delta_1\,.$$
Taking the $\inf$ over all $z\in \partial G_*$ in the latter
inequality, we obtain that $\widetilde{h}_*(x_0, \partial
G_*)\geqslant \frac{1}{2}\delta_1>0.$ This, $x_0\not\in \partial
G_*,$ as required. Now, $x_0\in G_*\subset D_*$ that contradicts
with $x_0\in
\partial D_*$ established above. The contradiction obtained above completes the proof.
\end{proof}

\section{Proof of Theorem~\ref{th1}}

The following statement may be found in~\cite[Lemma~4.3]{Vu}.

\begin{proposition}\label{pr1}
{\, Let $D$ be an open half space or an open ball in ${\Bbb C}$ and
let $E$ and $F$ be subsets of $D.$ Then
$$M_e(\Gamma(E, F, D))
\geqslant \frac{1}{2}\cdot M_e(\Gamma(E, F, \overline{{\Bbb
C}}))\,.$$ }
\end{proposition}

The following statement holds.

 \begin{lemma}\label{lem1}
{\, Let $a,$ $b,$ $c$ and $d$ be four different points of the
Riemannian surfaces of a domain $D$ in Riemannian surface ${\Bbb
S}.$ Then there are disjoint arcs $\gamma_1:[0,1]\rightarrow D$ and
$\gamma_2:[0, 1]\rightarrow D$ such that $\gamma_1(0)=a,$
$\gamma_1(1)=b,$ $\gamma_2(0)=c$ and $\gamma_2(1)=d.$}
 \end{lemma}

\begin{proof}
Let $p_0\in D$ and some normal neighborhood $U_0$ of $p_0$, such
that $\partial U_0\ne\varnothing.$ (Such neighborhoods exist due to
\cite[Proposition~18.1 and Lemma~18.3]{Sev}). Namely, we put
$U_0=\widetilde{B}(p_0, r_0).$ Let $\pi$ be a natural projection of
${\Bbb D}$ onto ${\Bbb S},$ and let $\varphi:=\pi^{\,-1}|_{U_0}.$
Now, $\varphi(\widetilde{B}(p_0, r_0))=B_h(x_0, r_0),$ where
$B_h(x_0, r_0)$ denotes, as usual, the hyperbolic ball in ${\Bbb
C},$ i.e., $B_h(x_0, r_0)=\{z\in {\Bbb C}: h(z, x_0)<r_0\}$ and $h$
is a hyperbolic metric, see~(\ref{eq7AC}). We may consider that
$x_0=0.$ Set $F:=\{a, b, c, d\}.$ By~\cite[Lemma~2.6,
Chapter~5]{Hir}, there exists a homeomorphism $\varphi_1$ of $D$
onto itself such that $\varphi_1(F)\subset \widetilde{B}(p_0,r_0).$
Observe that, the points $\varphi(\varphi_1(a))$ and
$\varphi(\varphi_1(b)),$ and $\varphi(\varphi_1(c))$ and
$\varphi(\varphi_1(d))$ may be joined by disjoint paths
$\widetilde{\gamma}_1$ and $\widetilde{\gamma}_2$ in $B_h(0, r_0),$
i.e., $\widetilde{\gamma}_1:[0,1]\rightarrow B_h(0, r_0)$ and
$\widetilde{\gamma}_2:[0, 1]\rightarrow B_h(0, r_0)$ such that
$\widetilde{\gamma}_1(0)=\varphi(\varphi_1(a)),$
$\widetilde{\gamma}_1(1)=\varphi(\varphi_1(b)),$
$\widetilde{\gamma}_2(0)=\varphi(\varphi_1(c))$ and
$\widetilde{\gamma}_2(1)=\varphi(\varphi_1(d))$ (see, e.g.,
\cite[Lemma~2.1]{SevSkv$_1$}). In this case, we set
$\gamma_1:=\varphi_1^{\,-1}(\pi(\widetilde{\gamma}_1))$ and
$\gamma_2:=\varphi_1^{\,-1}(\pi(\widetilde{\gamma}_2)).$ Observe
that, the paths $\gamma_1$ and $\gamma_2$ are desired.
\end{proof}

{\it Proof of Theorem~\ref{th1}.} \textbf{I.} Let us firstly prove
that the family $\frak{F}^{\delta}_{K, Q}(D, D_*)$ is equi-uniform
on $K$ under the integrability of $Q$ over concentric circles.
Assume the contrary. Now, there is $\varepsilon_0>0$ such that, for
any $m\in {\Bbb N}$ there are $x_m, y_m\in K$ and $f_m\in
\frak{F}^{\delta}_{K, Q}(D, D_*)$ such that $\widetilde{h}(x_m,
y_m)\geqslant \varepsilon_0,$ however, $\widetilde{h}_*(f_m(x_m),
f_m(y_m))<1/m.$ Since $K$ is a compactum, we may consider that
$x_m\rightarrow x_0$ and $y_m\rightarrow y_0$ as
$m\rightarrow\infty$ for some $x_0, y_0\in K.$ In addition, since
$\overline{{D}_*}$ is a compactum, we also may consider that
$f_m(x_m), f_m(y_m)\rightarrow w_0$ as $m\rightarrow\infty.$

\begin{figure}
\centerline{\includegraphics[scale=0.5]{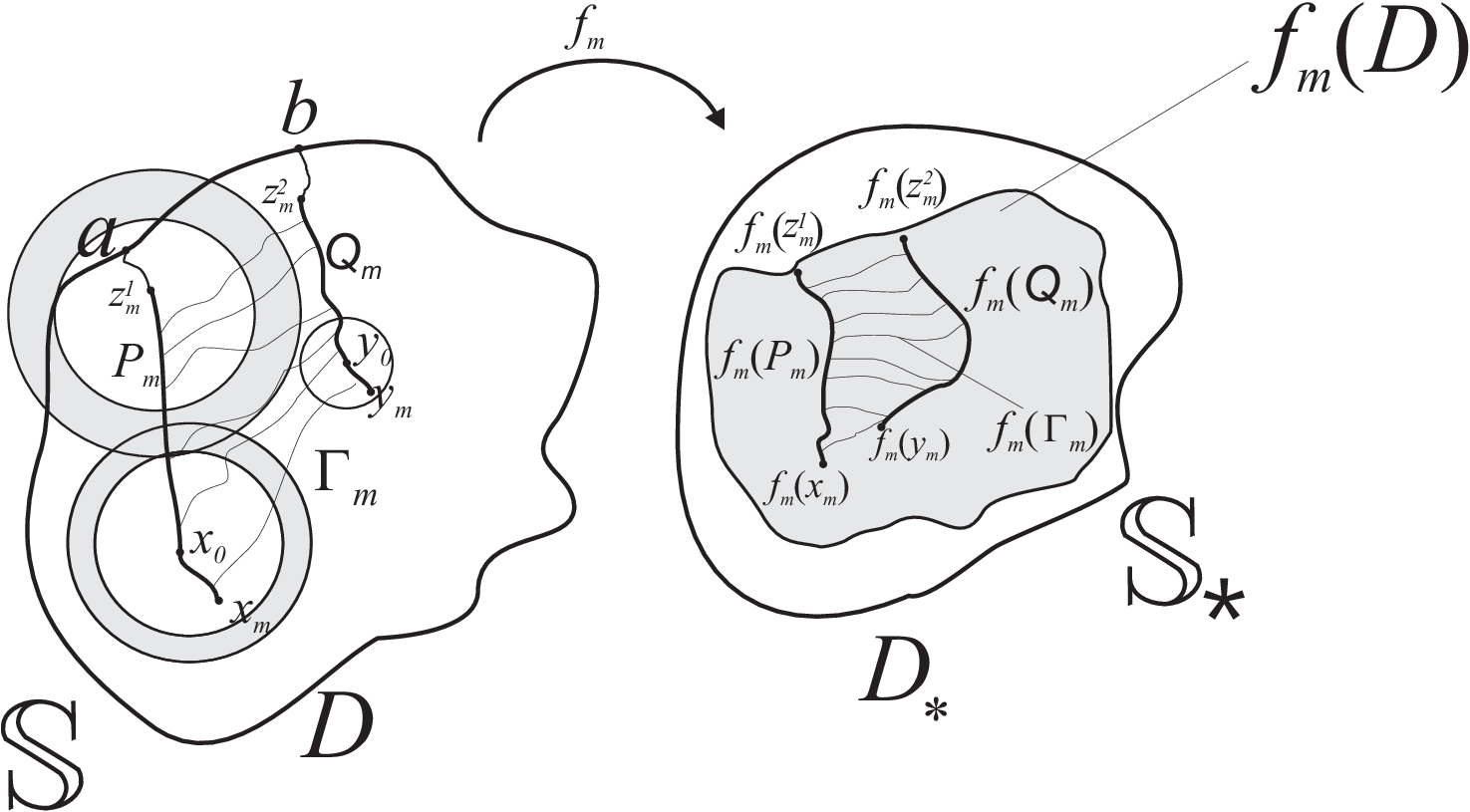}} \caption{To the
proof of Theorem~\ref{th1}, item~\textbf{I}}\label{figure2}
\end{figure}

\medskip
Let $a$ and $b$ be points from the condition of the theorem. By
Lemma~\ref{lem1} we may join points $x_0$ and $a$ and, respectively,
$y_0$ and $b$ by disjoint paths $\gamma_1:[1/2, 1]\rightarrow {\Bbb
S}$ and $\gamma_2:[1/2, 1] \rightarrow {\Bbb S}.$ Without loss of
generality, we may assume that both paths lie in $D,$ excepting $a$
and $b.$ Otherwise, due to~\cite[Theorem~1.I.5.46]{Ku$_2$}, we may
find subpaths of $\gamma_1$ and $\gamma_2$ of these paths with,
possibly, some other ends $\overline{x^*_1}$ and $\overline{x^*_2},$
but with the property mentioned above. So, we may consider below
that $\overline{x^*_1}=a$ and $\overline{x^*_2}=b.$ Let $R_1>0$ be
such that $\widetilde{B}(x_0, R_1)$ lies in some a normal
neighborhood of $x_0$ and $\overline{\widetilde{B}(x_0, R_1)}\cap
|\gamma_2|=\varnothing.$ In addition, let $R_2>0$ be such that
$\widetilde{B}(y_0, R_2)$ lies in the normal neighborhood of $y_0$
and
$$
(\overline{\widetilde{B}(x_0, R_1)}\cup|\gamma_1|)
\cap\overline{\widetilde{B}(y_0, R_2)}=\varnothing.
 $$
Since such balls $\widetilde{B}(x_0, R_1)$ and $\widetilde{B}(y_0,
R_2)$ lie in normal neighborhoods of points $x_0$ and $y_0,$
respectively, this balls are isometrically homeomorphic to the
corresponding hyperbolic disks $B_h(0, R_1)$ and $B_h(0, R_2)$ in
${\Bbb C}.$ Thus, $\widetilde{B}(x_0, R_1)$ and $\widetilde{B}(y_0,
R_2)$ are path connected. Moreover, $\widetilde{B}(x_0, r)$ and
$\widetilde{B}(y_0, r)$ are path connected for
$r\in[0,\max\{R_1,R_2\}].$ We may also assume that $x_m\in
\widetilde{B}(x_0, R_1)$ and $y_m\in \widetilde{B}(y_0, R_2)$ for
any $m\geqslant 1.$ Join the points $x_m$ and $x_0$ by a path
$\alpha^{\,*}_m:[0, 1/2] \rightarrow \widetilde{B}(x_0, R_1),$ and
join the point $y_m$ with the point $y_0$ by a path
$\beta^{\,*}_m:[0, 1/2]\rightarrow \widetilde{B}(y_0, R_2)$ (see
Figure~\ref{figure2}).

Set
$$
\alpha_m(t)=\quad\left\{\begin{array}{rr}
\alpha^*_m(t), & t\in [0, 1/2],\\
\gamma_1(t), & t\in [1/2, 1]\end{array} \right.\,,\quad
\beta_m(t)=\quad\left\{
\begin{array}{rr}
\beta^*_m(t), & t\in [0, 1/2],\\
\gamma_2(t), & t\in [1/2, 1].\end{array} \right.
 $$
By the construction, the sets
 $$
A_1:=|\gamma_1|\cup \overline{\widetilde{B}(x_0, R_1)}\,,\quad
A_2:=|\gamma_2|\cup \overline{\widetilde{B}(y_0, R_2)}
 $$
are disjoint so that there is $\varepsilon^{\,*}_1>0$ with
\begin{equation}\label{eq1J}
\widetilde{h}(A_1, A_2)=\inf\limits_{x\in A_1, y\in
A_2}\widetilde{h}(x, y)\geqslant \varepsilon^{\,*}_1>0\,.
 \end{equation}
 Let
$\delta(z)=\widetilde{h}(z, \partial U),$ where $U$ is such a normal
neighborhood of $z,$ and let
$$r_*(z):=\min\{\varepsilon^{\,*}_1, r_0(z), \delta(z)\}\,.$$
Cover the set $A_1$ with balls $\widetilde{B}(z, r_*/4),$ $z\in
A_1.$ Note that $|\gamma_1 |$ is a compact set in ${\Bbb S}$ as a
continuous image of the compactum $[1/2, 1]$ under the mapping
$\gamma_1.$ Then, by the Heine-Borel-Lebesgue lemma, there is a
finite subcover $\bigcup\limits_{i=1}^p\widetilde{B}(z_i, r_*/4)$ of
$A_1.$ Now,
\begin{equation}\label{eq2A}
A_1\subset \bigcup\limits_{i=1}^p\widetilde{B}(z_i,r_i/4)\,,\qquad
1\leqslant p<\infty\,,
 \end{equation}
where $r_i$ denotes the number $r_0(z_i)$ for any $z_i\in
\overline{D}.$ We may consider that $z_i\in D$ for all $1\leqslant
i\leqslant p$ because otherwise, if $z_i\not\in D,$ we may find
$z^{\,*}_i\in D$ and $r^{\,*}_i>0$ with $\widetilde{B}(z_i,
r_i/4)\subset\widetilde{B}(z^{\,*}_i,r^{\,*}_i/4).$

\medskip
Let $\alpha^0_m=f_m(\alpha_m)$ and $\beta^0_m=f_m(\beta_m).$ Since
$f_m$ is a homeomorphism, $\alpha^0_m(t_k)\rightarrow
\partial f_m(D)$ and $\beta^0_m(t^{\,\prime}_k)\rightarrow
\partial f_m(D)$ for any sequences $t_k\rightarrow 1-0$
and $t^{\,\prime}_k\rightarrow 1-0,$ because $C(f_m, \partial
D)\subset \partial f_m(D)$ (see, e.g.,
\cite[Proposition~13.5]{MRSY}). Now, there is $k=k_m\in {\Bbb N}$
such that $t_{k_m}, t^{\,\prime}_{k_m}\in (0, 1)$ such that
$\widetilde{h}_*(\alpha^0_m(t_{k_m}), \partial f_m(D))<\frac{1}{m},$
$\widetilde{h}_*(\beta^0_m(t^{\,\prime}_{k_m}), \partial
f_m(D))<\frac{1}{m}.$ Let $z^1_m=\alpha_m(t_{k_m}),$
$z^2_m=\beta_m(t^{\,\prime}_{k_m}),$ $P_m:=\bigl|\alpha_m|_{[0,
t_{k_m}]}\bigr|$ and $Q_m:=\bigl|\beta_m|_{[0,
t^{\,\prime}_{k_m}]}\bigr|.$ Now, by the triangle inequality and the
definition of the class $\frak{F}^{\delta}_{K, Q}(D, D_*),$
\begin{gather}\widetilde{h}_*(f_m(P_m))\geqslant
\widetilde{h}_*(f_m(\alpha_m(x_m)), \alpha^0_m(t_{k_m}))\nonumber\\
\geqslant \label{eq1}\widetilde{h}_*(f_m(\alpha_m(x_m)),
w)-\widetilde{h}_*(w, \alpha^0_m(t_{k_m}))\geqslant
\delta-\delta/2=\delta/2\,,
\end{gather}
while $w\in \partial f_m(D).$ Similarly,
\begin{equation}\label{eq2}\widetilde{h}_*(f_m(Q_m))\geqslant\delta/2\,,\qquad
m=1,2,\ldots\,. \end{equation}
Put $\Gamma_m:=\Gamma(P_m, Q_m, D).$ Recall that we write
$\Gamma_1>\Gamma_2$ if and only if each path $\gamma_1\in \Gamma_1$
has a subpath $\gamma_2\in \Gamma_2.$ (In other words, if $\gamma_1:
I\rightarrow{\Bbb S},$ then $\gamma_2: J\rightarrow{\Bbb S},$ where
$J\subset I$ and $\gamma_2(t)=\gamma_1(t)$ for $t\in J,$ and $I, J$
are segments, intervals, or half-intervals). Then, by~(\ref{eq1J})
and (\ref{eq2A}), and by~\cite[Theorem~1.I.5.46]{Ku$_2$}, we obtain
that
\begin{equation}\label{eq5}
\Gamma_m>\bigcup\limits_{i=1}^p\Gamma_{i}\,,
 \end{equation}
where $\Gamma_{i}:=\Gamma(\widetilde{S}(z_i, r_i/4),
\widetilde{S}(z_i, r_i/2), \widetilde{A}(z_i, r_i/4, r_i/2)).$ We
set $\widetilde{Q}(x)=\max\{Q(x), 1\}$ and
$$\widetilde{q}_{z_i}(r)=\int\limits_{\widetilde{S}(z_i,
r)}\widetilde{Q}(z)\,ds_{\widetilde{h}}(z)\,.$$ Now, we have also
that $\widetilde{q}_{z_i}(r)\ne \infty$ for any $r\in [r_i/4,
r_i/2].$
Set
$$I_i=I_i(z_i,r_i/4,r_i/2)=\int\limits_{r_i/4}^{r_i/2}\
\frac{dr}{r\widetilde{q}_{z_i}(r)}\,.$$
Observe that $I\ne 0,$ because $\widetilde{q}_{z_i}(r)\ne \infty$
for any $r\in E_1\subset [r_i/4,r_i/2].$ Besides that, note that
$I\ne\infty,$ since, for some constant $C>0,$
$$I_i\leqslant C\cdot\log\frac{\frac{r_i}{2}}{\frac{r_i}{4}}=
\log 2<\infty\,,\quad i=1,2, \ldots, p\,.$$
Now, we put
$$\eta_i(r)=\begin{cases}
\frac{1}{I_ir\widetilde{q}_{z_i}(r)}\,,&
r\in [r_i/4,r_i/2]\,,\\
0,& r\not\in [r_i/4,r_i/2]\,.
\end{cases}$$
Observe that, a function~$\eta_i$ satisfies the
condition~$\int\limits_{r_i/4}^{r_i/2}\eta_i(r)\,dr=1,$ therefore it
may be substituted in~(\ref{eq3*!!A}). In this case, we obtain that
\begin{equation}\label{eq7B}
M(f_m(\Gamma_{im}))\leqslant \int\limits_{\widetilde{A}(z_i, r_i/4,
r_i/2)} \widetilde{Q}(x)\,\eta^2_i(\widetilde{h}(x,
z_i))\,d\widetilde{v}(x)\,.\end{equation}
By Fubini's theorem for Riemannian surfaces (see, e.g.,
\cite[Lemma~17.1]{Sev}), we may calculate the expression in the
right-hand side of~(\ref{eq7B}). Indeed,
\begin{gather}\int\limits_{\widetilde{A}(z_i, r_i/4, r_i/2)}
\widetilde{Q}(x)\,\eta^2_i(\widetilde{h}(x,
z_i))\,d\widetilde{v}(x)\nonumber\\
\label{eq7C_1}
=\int\limits_{r_i/4}^{r_i/2}\int\limits_{\widetilde{S}(z_i,
r)}Q(x)\eta^2_i(\widetilde{h}(x,
z_i))\,ds_{\widetilde{h}}(x)\,dr\,\\
=\frac{1}{I_i^2}\int\limits_{r_i/4}^{r_i/2}r
\widetilde{q}_{z_i}(r)\cdot
\frac{dr}{r^2\widetilde{q}^{2}_{z_i}(r)}=\frac{1}{I_i}\,.\nonumber
\end{gather}
Now, by~(\ref{eq7B}) and~(\ref{eq7C_1}) we obtain that
$$M(f_m(\Gamma_{im}))\leqslant \frac{1}{I_i}\,,$$
whence from ~(\ref{eq5})
\begin{equation}\label{eq7D}
M(f_m(\Gamma_m))\leqslant
\sum\limits_{i=1}^pM(f_m(\Gamma_{im}))\leqslant
\sum\limits_{i=1}^p\frac{1}{I_i}:=C_0\,, \quad m=1,2,\ldots\,.
\end{equation}
On the other hand, we put $0<\delta_*<\delta/2$ such that
$\widetilde{B}(w_0, \delta_*)$ lies in some normal neighborhood of
$w_0.$ Now, by the definition, we may map $\widetilde{B}(w_0,
\delta_*)$ onto the hyperbolic ball $B_h(0, \delta_*)$ isometrically
by some isometry $\varphi.$ Due to~(\ref{eq1}) and (\ref{eq2}), and
by~\cite[Theorem~1.I.5.46]{Ku$_2$}, there are some subpaths
$\widetilde{\alpha}^0_m$ and $\widetilde{\beta}^0_m$ of $\alpha^0_m$
and $\beta^0_m$ starting at $f_m(x_m)$ and $f_m(y_m)$ which lie in
$\widetilde{B}(w_0, \delta_*)$ such that $\widetilde{\alpha}^0_m\cap
\widetilde{S}(w_0, \delta_*)\ne\varnothing\ne
\widetilde{\alpha}^0_m\cap \widetilde{S}(w_0, \delta_*)$ for
sufficiently large $m\in {\Bbb N}.$ Since all of the inner points of
$B_h(0, \delta_*)$ are weakly flat, see, e.g.,
\cite[Lemma~11.2]{Sev}, for any $P>0$ there is a number $m_0\in
{\Bbb N}$ such that
\begin{equation}\label{eq3}
M_e(\Gamma(\varphi(\widetilde{\alpha}^0_m),
\varphi(\widetilde{\beta}^0_m), B_h(0,
\delta_*)))>P\,,\end{equation}
where $M_e$ denotes the usual (Euclidean) modulus families of paths,
see \cite[Section~17.1]{Sev}. However, the Euclidean and the
hyperbolic moduli coincide, see~(\ref{eq10}). In other words,
\begin{gather}
M_e(\Gamma(\varphi(\widetilde{\alpha}^0_m),
\varphi(\widetilde{\beta}^0_m), B_h(0,
\delta_*)))=M_h(\Gamma(\varphi(\widetilde{\alpha}^0_m),
\varphi(\widetilde{\beta}^0_m), B_h(0,
\delta_*)))\nonumber\\\label{eq4}=M(\Gamma(\widetilde{\alpha}^0_m,
\widetilde{\beta}^0_m, \widetilde{B}(w_0, \delta_*)))\,.
\end{gather}
Finally, by the construction, $\Gamma(\widetilde{\alpha}^0_m,
\widetilde{\beta}^0_m, \widetilde{B}(w_0, \delta_*))\subset
f_m(\Gamma_m),$ so that
\begin{equation}\label{eq5A}
M(\Gamma(\widetilde{\alpha}^0_m, \widetilde{\beta}^0_m,
\widetilde{B}(w_0, \delta_*)))\leqslant M(f_m(\Gamma_m))\,.
\end{equation}
Combining~(\ref{eq3}), (\ref{eq4}) and~(\ref{eq5A}), we obtain that
\begin{equation}\label{eq6}
M(f_m(\Gamma_m))\geqslant P\,,\quad m\geqslant m_0\,.
\end{equation}
The relation~(\ref{eq6}) contradicts with~(\ref{eq7D}), which
finishes the proof of the equi-uniformness of the family
$\frak{A}^{\delta, r}_{A, Q}(D, D_*)$ that is the first part of
Theorem~\ref{th1}.

\medskip
\textbf{II.} We now prove the relation~(\ref{eq1D}) under the
additional conditions that (\ref{eq3*!!A})--(\ref{eq8BC}) hold for
any $p_0\in D,$ any $0<r_1<r_2<\sup\limits_{p\in D}\widetilde{h}(p,
p_0)$ and $Q\in L^1(D).$ It is convenient to carry out a proof by
contradiction, namely, assume that~(\ref{eq1D}) does not holds for
any $C_1.$ Now, given $k\in {\Bbb N}$ there is $x_k, y_k\in K,$
$x_k\ne y_k,$ and $f_k\in \frak{A}^{\delta, r}_{A, Q}(D, D_*)$ such
that
\begin{equation}\label{eq1E}
\widetilde{h}_*(f_k(x_k), f_k(y_k))<
\frac{1}{k}\cdot\exp\left\{-\frac{2\Vert
Q\Vert_1}{C{(\widetilde{h}(x_k, y_k))}^2}\right\}\,,
\end{equation}
where $C$ is a constant from~(\ref{eq3C}) for $n=2$ and $X={\Bbb
C}.$ Since $K$ is a compactum and $x_k, y_k\in K,$ we may consider
that $x_k\rightarrow x_0$ and $y_k\rightarrow y_0$ as $k\rightarrow
\infty,$ where $x_k, y_k\in K.$ We now consider two different cases.

\medskip
{\bf Case 1).} Let $x_0\ne y_0.$ Now, by the triangle inequality,
there is $\delta_0>0$ such that $\widetilde{h}(x_k,
y_k)\geqslant\delta_0$ for sufficiently large $k\in {\Bbb N}.$ In
this case, by the proving above (see item~\textbf{I}), there is
$\varepsilon_0>0$ such that $\widetilde{h}_*(f_k(x_k),
f_k(y_k))\geqslant \varepsilon_0$ for sufficiently large
$k=1,2,\ldots , $ as well. Again, by the triangle inequality,
$\widetilde{h}(x_k, y_k)\leqslant\delta_1$ for some $\delta_1>0$ and
sufficiently large $k\in {\Bbb N}.$ Now, the relation~(\ref{eq1E})
implies that
\begin{gather}
\varepsilon_0\leqslant\widetilde{h}_*(f_k(x_k), f_k(y_k))<
\frac{1}{k}\cdot\exp\left\{-\frac{2\Vert
Q\Vert_1}{C{(\widetilde{h}(x_k, y_k))}^2}\right\}\nonumber\\
\label{eq1F}\leqslant \frac{1}{k}\cdot\exp\left\{-\frac{2\Vert
Q\Vert_1}{C{\delta_1}^2}\right\}\,.
\end{gather}
It follows from~(\ref{eq1F}) that $\varepsilon_0<
\frac{1}{k}\cdot\exp\left\{-\frac{2\Vert
Q\Vert_1}{C{\delta_1}^2}\right\},$ or
$\varepsilon_0\cdot\exp\left\{\frac{2\Vert
Q\Vert_1}{C{\delta_1}^2}\right\}< \frac{1}{k}$ for sufficiently
large $k\in {\Bbb N},$ which is, obviously, impossible.

\begin{figure}
\centerline{\includegraphics[scale=0.4]{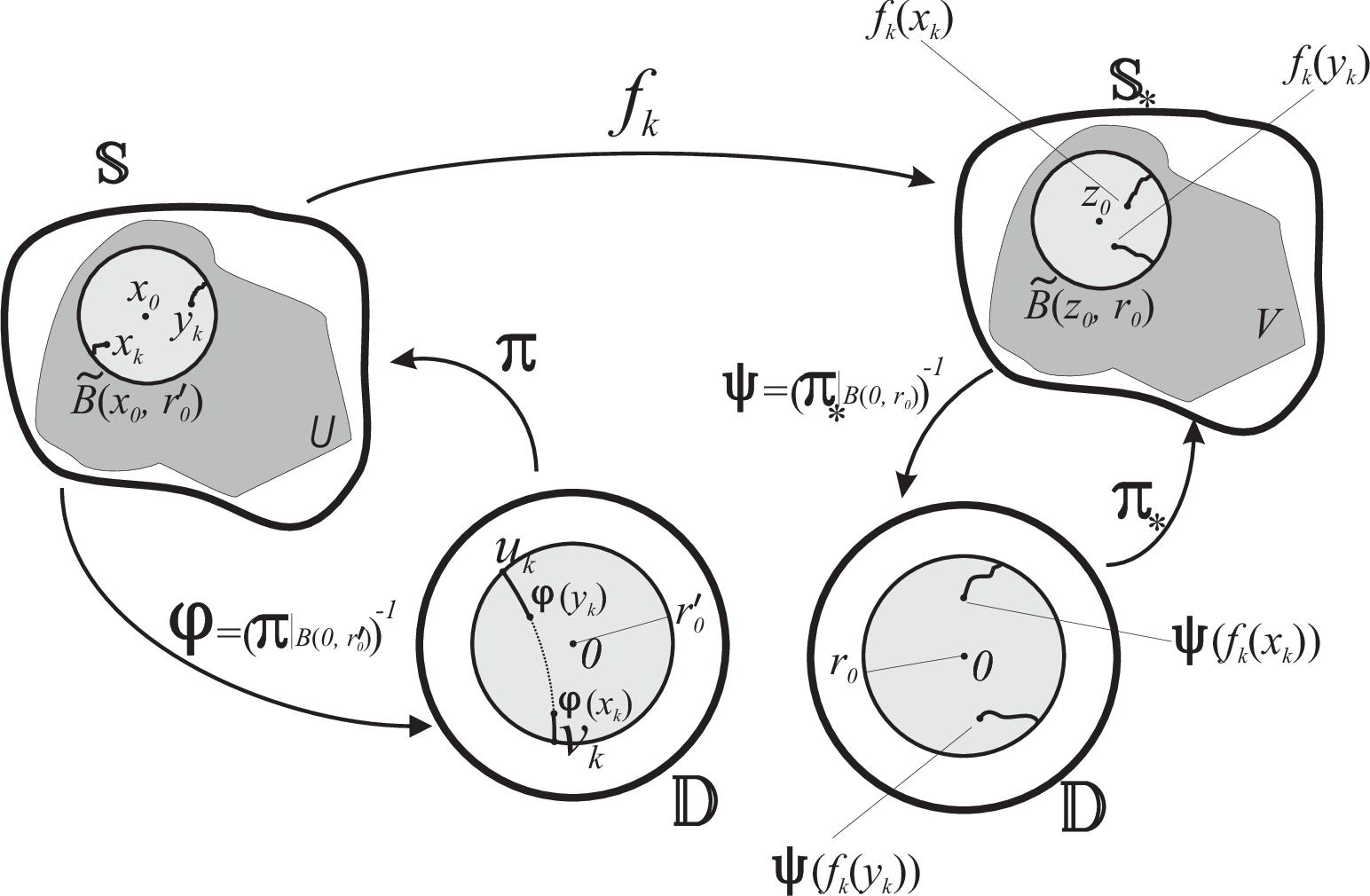}} \caption{To the
proof of Theorem~\ref{th1}, item~\textbf{II}, {\bf Case
2)}}\label{fig1}
\end{figure}

\medskip
{\bf Case 2).} Let $x_0=y_0.$ We choose some a normal neighborhood
$U$ of $x_0,$ see Figure~\ref{fig1}.
Similarly to Lemma~\ref{lem1}, let $\pi$ be a natural projection of
${\Bbb D}$ onto ${\Bbb S},$ and let $\varphi:=\pi^{\,-1}|_{U}.$ Let
$0<r^{\,\prime}_0<\widetilde{h}(x_0, \partial U).$ Now,
$\varphi(\widetilde{B}(x_0, r^{\,\prime}_0))=B_h(0,
r^{\,\prime}_0),$ where $B_h(0, r^{\,\prime}_0)=\{z\in {\Bbb C}:
h(z, 0)<r^{\,\prime}_0\}$ and $h$ is a hyperbolic metric. Since
$\overline{{D}_*}$ is a compactum, we also may consider that
$f_k(x_k)\rightarrow z_0$ and $f_k(y_k)\rightarrow z_1$ as
$k\rightarrow\infty.$ Observe that, $z_0=z_1,$ because otherwise it
follows from~(\ref{eq1E}) and by the triangle inequality that
\begin{equation}\label{eq1H}
\widetilde{h}_*(f_k(x_k), f_k(y_k))<
\frac{1}{k}\cdot\exp\left\{-\frac{2\Vert
Q\Vert_1}{C{(\widetilde{h}(x_k, y_k))}^2}\right\}\,,
\end{equation}
and the right side of~(\ref{eq1H}) converges to $0$ as
$k\rightarrow\infty;$ however, $\widetilde{h}_*(f_k(x_k),
f_k(y_k))\geqslant \frac{1}{2}\widetilde{h}_*(z_0, z_1)\ne 0$ for
sufficiently large $k=1,2,\ldots,$ that is impossible. Now,
$z_1=z_0,$ as required.

We choose some a normal neighborhood $V$ of $z_0.$ Let $\pi_*$ be a
natural projection of ${\Bbb D}$ onto ${\Bbb S}_*,$ and let
$\psi:=\pi_*^{\,-1}|_{V}.$ Let $0<r_0<\widetilde{h}_*(z_0,
\partial V).$ Now, $\psi(\widetilde{B}(z_0,
r_0))=B_h(0, r_0),$ where $B_h(0, r_0)=\{z\in {\Bbb C}: h(z,
0)<r_0\}$ and $h$ is a hyperbolic metric. We may consider that $r_0$
is sufficiently small such that $h(x, y)C_1(r_0)\leqslant |x-y|$ for
all $x, y\in B_h(0, r_0)$ and some constant $C_1=C_1(r_0)>0$ (see,
e.g., \cite[Lemma~18.4]{Sev}). Besides that, we may consider that
$\widetilde{B}(z_0, r_0)\subset f_k(D)$ for sufficiently large $k.$
Otherwise, $\widetilde{B}(z_0, 1/k)\cap f_k(D)\ne\varnothing\ne
\widetilde{B}(z_0, 1/k)\cap ({\Bbb S}_*\setminus f_k(D))$ for
sufficiently large $k\in {\Bbb N}.$
By~\cite[Theorem~1.I.5.46]{Ku$_2$}, $\widetilde{B}(z_0, 1/k)\cap
\partial f_k(D)\ne\varnothing,$ so that there is $w_k\in \partial
f_k(D)$ such that $\widetilde{h}_*(w_k, z_0)<1/k.$ However, by the
triangle inequality,
$$\widetilde{h}_*(w_k, f_k(x_k))\leqslant \widetilde{h}_*(w_k, z_0)
+\widetilde{h}_*(z_0, f_k(x_k))<1/k+\widetilde{h}_*(z_0,
f_k(x_k))\rightarrow 0$$
as $k\rightarrow\infty,$ i.e., $\widetilde{h}_*(f_k(K), \partial
f_k(D))\rightarrow 0,$ that contradicts the definition of the family
$\frak{F}^{\delta}_{K, Q}(D, D_*).$

\medskip
Let us draw a unique geodesic in the sense of the hyperbolic metric
through the points~$\varphi(x_k)$ and $\varphi(y_k)$ in ${\Bbb D}.$
Such a geodesic exists and is either an arc of the circle orthogonal
to the unit circle or the diameter of the unit disk in dependence on
the location of~$\varphi(x_k)$ and $\varphi(y_k)$; see, e.g.,
\cite[Section~7.3]{Be}. We parameterize this geodesic as
$\gamma_k:[0, 1]\rightarrow \overline{\Bbb D}$ such that
$\gamma_k(0)\in \partial {\Bbb D},$ $\gamma_k(1)\in \partial {\Bbb
D},$ besides that, $\gamma_k(t_1)=\varphi(x_k)$ and
$\gamma_k(t_2)=\varphi(y_k)$ for some $0<t_1<t_2<1.$ Observe that,
by the construction, $\bigl|\gamma_k|_{[0, t_1]}\bigr|\cap B_h(0,
r^{\,\prime}_0)\ne \varnothing\ne \bigl|{\gamma_k}_{[0,
t_1]}\bigr|\cap ({\Bbb D}\setminus B_h(0, r^{\,\prime}_0)).$ Thus,
by~\cite[Theorem~1.I.5.46]{Ku$_2$}, there is $0<p_1<t_1$ such that
$v_k:=\gamma_k(p_1)\in S_h(0, r^{\,\prime}_0)=\{z\in {\Bbb D}: h(0,
z)=r^{\,\prime}_0\}.$ Similarly, we may show that, there is
$t_2<p_2<1$ such that $u_k:=\gamma_k(p_2)\in S_h(0,
r^{\,\prime}_0).$ Let $\alpha_k:=\gamma_k|_{[p_1, t_1]}$ and
$\beta_k:=\gamma_k|_{[t_2, p_2]},$
$\widetilde{\alpha}_k=\pi(\alpha_k)$ and
$\widetilde{\beta}_k=\pi(\beta_k).$ Observe that,
$$h(v_k,
\varphi(y_k))=h(v_k, \varphi(x_k))+h(\varphi(x_k),
\varphi(y_k))\,,$$
see \cite[Theorem~7.3.2]{Be}. Now, since $\varphi$ is isometry,
\begin{equation*}\widetilde{h}(\pi(v_k), y_k)=\widetilde{h}(\pi(v_k),
x_k)+\widetilde{h}(x_k, y_k)\,.
\end{equation*}

{\bf III.} We now show that
\begin{equation}\label{eq1_1}
\widetilde{h}_*(f(\widetilde{S}(x_0, r^{\,\prime}_0)), \partial
f(D))\geqslant \delta_*>0
\end{equation}
for $f\in \frak{F}^{\delta}_{K, Q}(D, D_*).$

Otherwise, $\widetilde{h}_*(g_{k}(\widetilde{S}(x_0,
r^{\,\prime}_0)),
\partial g_{k}(D))<1/k,$ $k=1,2,\ldots ,$ $g_k\in \frak{F}^{\delta}_{K, Q}(D, D_*).$
Now, there is $b_k\in \widetilde{S}(x_0, r^{\,\prime}_0)$ such that
$\widetilde{h}_*(g_{k}(b_k),
\partial g_{k}(D))<1/k,$ $k=1,2,\ldots .$ Let $R_0>r_0^{\,\prime}>0$ be a number
for which $G:=\widetilde{B}(x_0, R_0)$ lies in the normal
neighborhood $U.$ We join $x_k$ and $b_k$ by a path $\alpha_k$
inside $\overline{\widetilde{B}(x_0, r^{\,\prime}_0)}$ (it is
possible because $\overline{\widetilde{B}(x_0, r^{\,\prime}_0)})$ is
some hyperbolic disk up to isometry). Let $A_k:=|\alpha_k|.$ Now,
since $x_k\in K,$ we have that $\widetilde{h_*}(g_k(x_k), \partial
g_{k}(D))\geqslant \delta.$ In addition, since
$\widetilde{h}_*(g_{k}(b_k),
\partial g_{k}(D))=\widetilde{h}_*(g_{k}(b_k),
c_k)<1/k$ for some $c_k\in \partial g_{k}(D),$ $k=1,2,\ldots ,$ we
have $\widetilde{h}_*(g_{k}(b_k), c_k)<\delta/2$ for sufficiently
large $k=1,2,\ldots .$ Now, by the triangle inequality, for
sufficiently large $k=1,2,\ldots ,$
\begin{gather}\widetilde{h}_*(g_k(A_k))\geqslant \widetilde{h}_*(g_k(x_k),
g_k(b_k))\nonumber\\
\label{eq1_3} \geqslant \widetilde{h}_*(g_k(x_k),
c_k)-\widetilde{h}_*(g_{k}(b_k), c_k)\geqslant
\delta-\delta/2=\delta/2\,.
\end{gather}
We show that
\begin{equation}\label{eq1_2}
\widetilde{h}_*(g_k(\widetilde{S}(x_0,
R_0)))=\widetilde{h}_*(\partial(g_k(\widetilde{B}(x_0,
R_0))))\geqslant \delta_{**}>0
\end{equation}
for some $\delta_{**}>0$ and sufficiently large $k=1,2,\ldots. $
Otherwise, if~(\ref{eq1_2}) is not true, we may consider that, there
is $z_*\in \overline{D_*}$ such that, for every $\sigma>0$ there is
$k=k(\sigma)\in {\Bbb N}$ with
\begin{equation}\label{eq1_4}
g_k(\widetilde{S}(x_0, R_0))\subset \widetilde{B}(z_*, \sigma)\,.
\end{equation}
Let $0<\sigma<\delta/2.$ Now, due to~(\ref{eq1_3}), it follows that
there is $\omega_k\in g_k(\widetilde{B}(x_0, R_0))\setminus
\widetilde{B}(z_*, \sigma).$ It follows from the assumptions of
Theorem~\ref{th1} that $\partial g_k(D)\ne\varnothing.$ Let
$\Delta_k\in \partial g_k(D).$ Now, we join the points $\omega_k$
and $\Delta_k$ by a path $\beta_k$ in ${\Bbb S}_*\setminus
\widetilde{B}(z_*, \sigma).$ (The latter is possible due to the path
connectedness of  $\widetilde{S}(z_*, \sigma)$ for sufficiently
small $\sigma>0.$) Since $|\beta_k|\cap g_k(\widetilde{B}(x_0,
R_0))\ne \varnothing \ne |\beta_k|\setminus g_k(\widetilde{B}(x_0,
R_0)),$ by \cite[Theorem~1.I.5.46]{Ku$_2$} $|\beta_k|\cap \partial
g_k(\widetilde{B}(x_0, R_0))\ne \varnothing.$ Since $\partial
g_k(\widetilde{B}(x_0, R_0))=g_k(\widetilde{S}(x_0, R_0))$ and
$\beta_k$ lies outside $ \widetilde{B}(z_*, \sigma),$ the latter
contradicts with~(\ref{eq1_4}). This contradiction
proves~(\ref{eq1_2}).

\medskip
Finally, we apply Lemma~\ref{lem2} to the domain
$D:=\widetilde{B}(x_0, R_0)$ and a continuum
$A:=\overline{\widetilde{B}(x_0, r^{\,\prime}_0)}.$ By this lemma,
$\widetilde{h}_*(g_k(A),
\partial g_k(G))\geqslant \delta_{***}$ for any $k=1,2,\ldots $ and some $\delta_{***}>0.$
By Proposition~\ref{pr6} $$\widetilde{h}_*(g_k(A),
\partial g_k(D))\geqslant \delta_{***}$$ for any $k=1,2,\ldots $ and some $\delta_{****}>0.$
The latter contradicts with our assumption that
$$\widetilde{h}_*(g_{k}(\widetilde{S}(x_0, r^{\,\prime}_0)),
\partial g_{k}(D))<1/k\,,\qquad k=1,2,\ldots .$$ This proves~(\ref{eq1_1}).

\medskip
{\bf IV.} By the proving above in~(\ref{eq1_1}),
$\frak{F}^{\delta}_{K, Q}(D, D_*)\subset
\frak{F}^{\widetilde{\delta}}_{\widetilde{K}, Q}(D, D_*)$ for,
possibly, some another $\widetilde{\delta}>0$ and
$\widetilde{K}=K\cup \widetilde{S}(x_0, r^{\,\prime}_0).$ Now,
observe that, by item~\textbf{I}, $\widetilde{h}_*(f_k(u_k),
z_0)\geqslant \delta_1>0$ and $\widetilde{h}_*(f_k(v_k),
z_0)\geqslant \delta_1>0$ for some $\delta_1>0$ and sufficiently
large $k\in{\Bbb N}.$ If $f_k(\widetilde{\alpha}_k)\subset
\widetilde{B}(z_0, r_0),$ we put $P_k:=f_k(\widetilde{\alpha}_k).$
Otherwise, $f_k(\widetilde{\alpha}_k)\cap \widetilde{S}(z_0,
r_0)\ne\varnothing,$ see \cite[Theorem~1.I.5.46]{Ku$_2$}. In this
case, we may find a subpath $P_k$ of $f_k(\widetilde{\alpha}_k)$
starting at $f_k(x_k)$ such that $\widetilde{h}_*(P_k)\geqslant
r_0/2.$ In any case, there is a subpath $P_k$ of
$f_k(\widetilde{\alpha}_k)$ starting at $f_k(x_k)$ with
$\widetilde{h}_*(P_k)\geqslant \min\{\delta_1, r_0/2\}:=\delta_0.$
Similarly, we may find a subpath $Q_k$ of $f_k(\widetilde{\beta}_k)$
starting at $f_k(y_k)$ such that $\widetilde{h}_*(Q_k)\geqslant
\delta_0.$ Now, paths $\psi(P_k)$ and $\psi(Q_k)$ are starting at
$\psi(f_k(x_k))$ and $\psi(f_k(y_k)),$ respectively, besides that,
$\psi(P_k)$ and $\psi(Q_k)$ lie in $B_h(0, r_0),$ and
$h(\psi(P_k))\geqslant \delta_0$ and $h(\psi(Q_k))\geqslant
\delta_0$ for sufficiently large $k=1,2,\ldots .$ Now,
$d(\psi(P_k))\geqslant \delta^{\,*}_0$ and $d(\psi(Q_k))\geqslant
\delta^{\,*}_0$ for sufficiently large $k=1,2,\ldots ,$ where $d(A)$
denotes the Euclidean diameter of the set $A\subset {\Bbb C}.$

\medskip
Now, we use the fact that ${\Bbb C}$ is a Loewner space (see, e.g.,
\cite[Theorem~8.2]{He}). Observe that,
%
%
$$\Delta(\psi(P_k), \psi(Q_k)):=\frac{d(\psi(P_k), \psi(Q_k))}
{\min\{d(\psi(P_k)), d(\psi(Q_k))\}} \leqslant
\frac{|\psi(f_k(x_k))-\psi(f_k(y_k))|}{\delta^{\,*}_0}$$
%
because $d(\psi(P_k), \psi(Q_k))=\inf\limits_{x\in \psi(P_k), y\in
\psi(Q_k)}|x-y|\leqslant |\psi(f_k(x_k))-\psi(f_k(y_k))|.$
Then, by the definition of the Loewner function $\Phi_n(t),$ $n=2,$
in~(\ref{eq2H}), we obtain that
\begin{equation}\label{eq4C}
\Phi_2\left(\frac{|\psi(f_k(x_k))-\psi(f_k(y_k))|}{\delta^{\,*}_0}\right)\leqslant
M_e(\Gamma(\psi(P_k), \psi(Q_k), {\Bbb C}))\,.
\end{equation}
Since $f_k(x_k), f_k(y_k)\rightarrow z_0$ as $k\rightarrow\infty,$
$\frac{|\psi(f_k(x_k))-\psi(f_k(y_k))|}{\delta^{\,*}_0}<\widetilde{\delta}_0,$
where $\widetilde{\delta}_0$ corresponds to~(\ref{eq3C}). Now,
by~(\ref{eq3C})
\begin{equation}\label{eq5C}
C\cdot\log
\frac{\delta^{\,*}_0}{|\psi(f_k(x_k))-\psi(f_k(y_k))|}\leqslant
\Phi_2\left(\frac{|\psi(f_k(x_k))-\psi(f_k(y_k))|}{\delta^{\,*}_0}\right)\,.
\end{equation}
Now, by~(\ref{eq4C}) and~(\ref{eq5C}) we obtain that
\begin{equation}\label{eq6C}
C\cdot\log
\frac{\delta^{\,*}_0}{|\psi(f_k(x_k))-\psi(f_k(y_k))|}\leqslant
M_e(\Gamma(\psi(P_k), \psi(Q_k), {\Bbb C}))\,.
\end{equation}
On the other hand, by Proposition~\ref{pr1}, since $B_h(0, r_0)$ is
some Euclidean ball centered at the origin, it follows
from~(\ref{eq6C}) that
\begin{gather}
M_e(\Gamma(\psi(P_k), \psi(Q_k), {\Bbb C}))\leqslant 2
M_e(\Gamma(\psi(P_k), \psi(Q_k), B_h(0, r_0)))\nonumber\\
\label{eq6D_1}= 2M_h(\Gamma(\psi(P_k), \psi(Q_k), B_h(0, r_0)))=
2M(\Gamma(P_k, Q_k, \widetilde{B}(0, r_0)))\\
\nonumber\leqslant 2M(\Gamma(P_k, Q_k, f_k(D)))\,.
\end{gather}
Here we took into account that the Euclidean and the hyperbolic
moduli coincide, see, e.g., \cite[relation~(9.31)]{Sev}. Since
$$|\psi(f_k(x_k))-\psi(f_k(y_k))|\leqslant
h(\psi(f_k(x_k)),\psi(f_k(y_k)))=\widetilde{h}_*(f_k(x_k),
f_k(y_k)),$$ from~(\ref{eq6C}) and (\ref{eq6D_1}) we obtain that
\begin{equation}\label{eq1B}
\frac{C}{2}\log \frac{\delta^{\,*}_0}{\widetilde{h}_*(f_k(x_k),
f_k(y_k))}\leqslant M(\Gamma(P_k, Q_k, f_k(D)))\,.
\end{equation}
Now let us prove the upper bound for $M(\Gamma(P_k, Q_k, f_k(D))).$
Indeed, by the definition,
\begin{equation}\label{eq2_1}
\Gamma(P_k, Q_k, f_k(D))=f_k(\Gamma(P^{\,*}_k, Q^{\,*}_k, D))\,,
\end{equation}
where $P^{\,*}_k:=f_k^{\,-1}(P_k)$ and $Q^{\,*}_k:=f_k^{\,-1}(Q_k).$
In addition, we now show that
\begin{equation}\label{eq1A}
\Gamma(P^{\,*}_k, Q^{\,*}_k, D)>\Gamma(\widetilde{S}(\pi(v_k),
\varepsilon^{(1)}), \widetilde{S}(\pi(v_k), \varepsilon^{(2)}),
\widetilde{A}(\pi(v_k), \varepsilon^{(1)}, \varepsilon^{(2)}))\,,
\end{equation}
where $\varepsilon^{(1)}:=\widetilde{h}(\pi(v_k), x_k)$ and
$\varepsilon^{(2)}:=\widetilde{h}(\pi(v_k), y_k)$ Observe that, by
the simplest properties of geodesics (see, e.g.,
\cite[Theorems~7.3.1 and 7.3.2]{Be} )
\begin{gather}
\widetilde{h}(x_k, y_k)+\varepsilon^{(1)}\nonumber\\
\label{eq5B}=\widetilde{h}(x_k, y_k)+\widetilde{h}(x_k, \pi(v_k))=
\widetilde{h}(y_k, \pi(v_k))=\varepsilon^{(2)}\,,
\end{gather}
and, thus, $\varepsilon^{(1)}<\varepsilon^{(2)}.$ Let
$\gamma\in\Gamma(P^{\,*}_k, Q^{\,*}_k, D).$ Then $\gamma:[0,
1]\rightarrow D,$ $\gamma(0)\in P^{\,*}_k,$ $\gamma(1)\in Q^{\,*}_k$
and $\gamma(s)\in D$ for $0<s<1.$ Since $\gamma(0)\in P^{\,*}_k,$
$\gamma(0)\in \widetilde{B}(\pi(v_k), \varepsilon^{(1)}).$ Since
$\gamma(1)\in Q^{\,*}_k,$ $\gamma(1)\in {\Bbb S}\setminus
\widetilde{B}(\pi(v_k), \varepsilon^{(2)}).$ Now,
by~\cite[Theorem~1.I.5.46]{Ku$_2$} there is $\zeta_1\in (0, 1)$ such
that~$\gamma(\zeta_1)\in \widetilde{S}(\pi(v_k),
\varepsilon^{(1)}).$ We may assume that $\gamma(t)\not\in
\widetilde{B}(\pi(v_k), \varepsilon^{(1)})$ for $t\in (\zeta_1, 1).$
Put $\alpha^1:=\gamma|_{[\zeta_1, 1]}.$ On the other hand, since
$\varepsilon^{(1)}<\varepsilon^{(2)}$ and $\gamma(\zeta_1)\in
\widetilde{S}(\pi(v_k), \varepsilon^{(1)}),$ we obtain that
$|\alpha^1|\cap B(\pi(v_k), \varepsilon^{(2)})\ne\varnothing.$
Again, by~\cite[Theorem~1.I.5.46]{Ku$_2$} there is $\zeta_2\in
(\zeta_1, 1)$ such that~$\gamma(\zeta_2)\in \widetilde{S}(\pi(v_k),
\varepsilon^{(2)}).$ Set $\alpha^2:=\alpha^1|_{[\zeta_1, \zeta_2]}.$
Now, $\gamma>\alpha^2$ and $\alpha^2\in \Gamma(S(\pi(v_k),
\varepsilon^{(1)}), \widetilde{S}(\pi(v_k), \varepsilon^{(2)}),
\widetilde{A}),$ $\widetilde{A}=\widetilde{A}(\pi(v_k),
\varepsilon^{(1)}, \varepsilon^{(2)}).$ Thus, (\ref{eq1A}) is
proved.

\medskip
By the inequalities~(\ref{eq1B}), (\ref{eq2_1}), (\ref{eq1A}), we
obtain that
\begin{gather}
\frac{C}{2}\cdot\log \frac{\delta^{\,*}_0}{\widetilde{h}_*(f_k(x_k),
f_k(y_k))}\nonumber\\ \label{eq15}\leqslant
M(f(\Gamma(\widetilde{S}(\pi(v_k), \varepsilon^{(1)}),
\widetilde{S}(\pi(v_k), \varepsilon^{(2)}), \widetilde{A}(\pi(v_k),
\varepsilon^{(1)}, \varepsilon^{(2)}))))\,.
\end{gather}
Now, we set
$$\eta(t)= \left\{
\begin{array}{rr}
\frac{1}{\widetilde{h}(x_k, y_k)}, & t\in [\varepsilon^{(1)}, \varepsilon^{(2)}],\\
0,  &  t\not\in [\varepsilon^{(1)}, \varepsilon^{(2)}]\,.
\end{array}
\right.$$
Observe that, $\eta$ satisfies the relation~(\ref{eq8BC}) for
$r_1=\varepsilon^{(1)}$ and $r_2=\varepsilon^{(2)}.$ Indeed, it
follows from~(\ref{eq5B}) that
$$r_2-r_1=\varepsilon^{(2)}-\varepsilon^{(1)}=\widetilde{h}(x_k, y_k)\,.$$
Then
$\int\limits_{\varepsilon^{(1)}}^{\varepsilon^{(2)}}\eta(t)\,dt=(1/\widetilde{h}(x_k,
y_k))\cdot (\varepsilon^{(2)}-\varepsilon^{(1)})= 1.$
By~(\ref{eq15}) and by the definition of the class
$\frak{F}^{\delta}_{K, Q}(D, D_*),$ we obtain that
\begin{equation}\label{eq14***A}
\frac{C}{2}\cdot\log \frac{\delta^{\,*}_0}{\widetilde{h}_*(f_k(x_k),
f_k(y_k))}\leqslant \frac{1}{{(\widetilde{h}(x_k,
y_k))}^2}\int\limits_{D} Q(w)\,d\widetilde{v}(w)=\frac{\Vert
Q\Vert_1}{{(\widetilde{h}(x_k, y_k))}^2}\,.
\end{equation}
Finally, (\ref{eq14***A}) implies that
\begin{equation}\label{eq3A}
\widetilde{h}_*(f_k(x_k), f_k(y_k))\geqslant
\delta^{\,*}_0\exp\left\{-\frac{2\Vert
Q\Vert_1}{C{(\widetilde{h}(x_k, y_k))}^2}\right\}\,,
\end{equation}
that contradicts with~(\ref{eq1E}). The theorem is completely
proved.~$\Box$

\medskip
{\it Proof of Corollary~\ref{cor1}} follows by the Fubini theorem,
see, e.g., \cite[Lemma~17.1]{Sev}. Indeed, by the above lemma, for
sufficiently small $r_0>0,$
$$
\int\limits_{\widetilde{B}(p_0,
r_0)}Q(p)\,d\widetilde{v}(p)=\int\limits_0^{r_0}\int\limits_{\widetilde{S}(p_0,
r)}Q(p)\,ds_{\widetilde{h}}(p)\,dr\,,$$
whence the existence and the finiteness of the integral
$\int\limits_{\widetilde{S}(p_0, r)}Q(p)\,ds_{\widetilde{h}}(p)$
follows for almost all $r\in (0, r_0)$ follows.~$\Box$

\medskip
{\it Proof of Corollary~\ref{cor2}.} We will argue similarly to the
case ${\Bbb R}^n,$ see, e.g., \cite[Theorem~1.3]{RS$_1$}. Assume the
contrary, i.e., the family $\frak{F}^{\delta}_{K, Q}(D, D_*)$ is not
uniformly open on $K.$ Then there exists $\varepsilon_0>0$ such
that, for any $m\in {\Bbb N}$ there exists $x_m\in K$ and $f_m\in
\frak{F}^{\delta}_{K, Q}(D, D_*)$ with $\widetilde{B}(x_m,
\varepsilon_0)\subset K$ while $\widetilde{B}(f_m(x_m),
1/m)\setminus f_m(\widetilde{B}(x_m, \varepsilon_0))\ne\varnothing.$
Let $y_m\in \widetilde{B}(f_m(x_m), 1/m)\setminus f_m(B(x_m,
\varepsilon_0)).$ Since by the assumption $\overline{D}_*$ is a
compactum in ${\Bbb S}_*,$ we may consider that $f_m(x_m)$ and $y_m$
converge to some point $\omega_{*}\in\overline{D_*}$ as
$m\rightarrow \infty.$ Now, $f_m(x_m)$ and $y_m$ lie in the disk
which is isometrically equivalent to the hyperbolic disk on the
plane. Arguing similarly to the proof of Theorem~\ref{th1}, we may
join the points $f_m(x_m)$ and $y_m$ by the geodesic line
$r_m=r_m(t),$ $t\in [0, 1],$ such that $r_m(0)=f_m(x_m)$ and
$r_m(1)=y_m.$ Since $|r_m|\cap f_m(\widetilde{B}(x_m,
\varepsilon_0))\ne\varnothing\ne |r_m|\setminus
f_m(\widetilde{B}(x_m, \varepsilon_0)),$
by~\cite[Theorem~1.I.5.46]{Ku$_2$} there is a point $z_m=r_m(t_m)\in
\partial f_m(\widetilde{B}(x_m, \varepsilon_0)).$ Without loss of generality, we
may assume that the path $\beta_m:=r_m|_{[0, t_m)}$ lies in
$f_m(\widetilde{B}(x_m, \varepsilon_0)).$ Put
$\alpha_m=f_m^{\,-1}(\beta_m).$ Since homeomorphisms preserve the
boundary (see, e.g., Proposition~13.5 in \cite{MRSY}),
$\alpha_m(t)\rightarrow \widetilde{S}(x_m, \varepsilon_0)$ as
$t\rightarrow t_m.$ Now, there is a point $z_m\in |\alpha_m|$ with
$\widetilde{h}(z_m, x_m)>\varepsilon_0/2.$ Let
$\widetilde{\alpha}_m$ be a part of $\alpha_m$ starting at $x_m$ and
ending at $z_m.$ We set $E_m:=|\widetilde{\alpha}_m|.$

\medskip
Finally, $\widetilde{h}(E_m)\geqslant \varepsilon_0/2,$ $E_m\subset
K$ for any $m=1,2,\ldots $ and, simultaneously,
$\widetilde{h}_*(f_m(E_m))\leqslant \widetilde{h}_*(|\beta_m|)<1/m$
for sufficiently large $m=1,2,\ldots ,$ for
$f_m\in\frak{F}^{\delta}_{K, Q}(D, D_*).$ The latter contradicts the
statement of Theorem~\ref{th1}. The obtained contradiction completes
the proof of the theorem.~$\Box$

\begin{theorem}\label{th1A}
{\, Let $D$ and $D_*$ be domains in ${\Bbb S}$ and ${\Bbb S}_*,$
respectively, while $D$ has at least two different boundary points
$a$ and $b$ and, besides that, $\overline{D_*}$ is a compactum in
${\Bbb S}_*.$ Assume that, for each point $p_0\in D$ there is
$r_0=r_0(p_0)>0,$ $r_0<\sup\limits_{p\in D}\widetilde{h}(p, p_0),$
such that for any $0<r_1<r_2<r_0$ there is a set $E_1\subset[r_1,
r_2]$ of a positive linear Lebesgue measure such that the function
$Q$ is integrable over $\widetilde{S}(p_0, r)$ for every $r\in E_1.$
Let $K$ be a compactum in $D.$ Then the family $\frak{A}^{\delta,
r}_{A, Q}(D, D_*)$ is uniformly light on $K;$ moreover,
$\frak{A}^{\delta, r}_{A, Q}(D, D_*)$ is equi-uniform on $K.$

\medskip
In addition, if (\ref{eq3*!!A})--(\ref{eq8BC}) hold for any $p_0\in
D,$ any $0<r_1<r_2<\sup\limits_{p\in D}\widetilde{h}(p, p_0),$ any
$f\in \frak{A}^{\delta, r}_{A, Q}(D, D_*)$ and $Q\in L^1(D),$ then
there exists $C_1>0$ such that
\begin{equation*}
\widetilde{h}_*(f(x), f(y))\geqslant
C_1\cdot\exp\left\{-\frac{2\Vert Q\Vert_1}{C{(\widetilde{h}(x,
y))}^2}\right\}
\end{equation*}
for all $x, y\in K,$ $x\ne y,$ and every $f\in \frak{A}^{\delta,
r}_{A, Q}(D, D_*),$ where $C$ is a constant from (\ref{eq3C}) for
$n=2$ and $X={\Bbb C}.$ }
\end{theorem}

\begin{proof} Let $K$ be a compactum in $D.$ By Lemma~\ref{lem2} and
Proposition~\ref{pr4}, there is $\delta^{\,*}_1>0$ such that
$h(f(K),
\partial f(D))\geqslant \delta^{\,*}_1$ for any $f\in\frak{A}^{\delta, r}_{A, Q}(D).$
Thus, $f\in\frak{F}^{\delta^{\,*}_1}_{K, Q}(D, D_*).$ In this case,
the desired statement follows directly from Theorem~\ref{th1}.
\end{proof}

\begin{corollary}\label{cor3}
{\, The statement of Theorem~\ref{th1A} remains true if we replace
the assumption on $Q$ mentioned above concerning the integrability
of $Q$ over circles $\widetilde{S}(p_0, r)$ by a simpler condition:
$Q\in L_{\rm loc}^1(D).$}
\end{corollary}

The proofs of the following Corollaries are completely similar to
the proofs of Corollary~\ref{cor1} and~\ref{cor2}, respectively.

\begin{corollary}\label{cor5}{\bf \,(An analogue of Koebe-Bloch
theorem)}.~{ Under assumptions of Theorem~\ref{th1A}, the family
$\frak{A}^{\delta, r}_{A, Q}(D, D_*)$ is uniformly open on every
compactum $K$ in $D.$ In particular, the latter is true if we
replace the assumption on $Q$ mentioned in Theorem~\ref{th1A}
concerning the integrability of $Q$ over circles $\widetilde{S}(p_0,
r)$ by a simpler condition: $Q\in L_{\rm loc}^1(D).$}
\end{corollary}

\section{Another lemma on approaching continua}

Firstly we prove some respectively easy statement.

\begin{proposition}\label{pr2}
Let $A$ be a compactum in ${\Bbb S}.$ Given $x_0\in A$ we denote by
$r_*(x_0)$ the supremum over all $r>0$ for which $\widetilde{B}(x_0,
r)$ lies in the normal neighborhood $U$ of $x_0.$ Now, the following
is true: $r_*:=\inf\limits_{x_0\in A}r_*(x_0)>0.$
\end{proposition}

\begin{proof}
Proving this statement from the contradiction, we assume that there
is a sequence $x_m\in A$ such that $r_*(x_m)<1/m,$ $m=1,2,\ldots .$
Since $A$ is a compactum in ${\Bbb S},$ we may consider that
$x_m\rightarrow x_0$ as $m\rightarrow\infty.$ Let $U$ ne a normal
neighborhood of $x_0.$ Now, $x_m\in U$ for sufficiently large $m\in
{\Bbb N}.$ Obviously, $r_*(x_m)>r_*(x_0)/2$ for such $m,$ because
$\widetilde{B}(x_m, r_*(x_0)/2)\subset \widetilde{B}(x_0, r_*(x_0))$
if $m$ is sufficiently large. The above contradiction finishes the
proof.
\end{proof}

\medskip
The following assertion is true.

\begin{lemma}\label{lem4}{\, Let $D$ and $D_*$ be domains in ${\Bbb S}$ and ${\Bbb S}_*,$
respectively, while $\overline{D_*}$ is a compactum in ${\Bbb S}_*$.
Let $f_m\in \frak{B}^{\delta}_{A, Q}(D, D_*).$ Assume that, for each
point $p_0\in D$ there is $r_0=r_0(p_0)>0,$ $r_0<\sup\limits_{p\in
D}\widetilde{h}(p, p_0),$ such that for any $0<r_1<r_2<r_0$ there is
a set $E_1\subset[r_1, r_2]$ of a positive linear Lebesgue measure
such that the function $Q$ is integrable over $\widetilde{S}(p_0,
r)$ for every $r\in E_1.$ If $\partial f_m(D)\ne \varnothing,$ then
there exists $\delta_1>0$ such that $\widetilde{h}_*(f_m(A),
\partial f_m(D))\geqslant \delta_1$ for any $m=1,2,\ldots .$ }
\end{lemma}

\begin{proof} We argue similarly to Lemma~3.3 in \cite{RS$_1$}.
By proving by contradiction, i.e., assuming that the desired
statement is not true, we assume that there exists an increasing
subsequence of numbers $m_k,$ $k=1,2,\ldots ,$ such that
$\widetilde{h}_*(f_{m_k}(A),
\partial f_{m_k}(D))\rightarrow 0$ as $k\rightarrow\infty.$ Without loss of generality, we
may assume that $f_m$ satisfy the above condition, i.e.
\begin{equation}\label{eq1_9}
\widetilde{h}_*(f_m(A),
\partial f_m(D))\rightarrow 0\,,\qquad m\rightarrow\infty\,.
\end{equation}
Otherwise, we may consider the sequence $f_{m_k}$ instead of $f_m$
in the conditions of the lemma.

\medskip
Let $x_m, y_m$ be points in $A$ such that $\widetilde{h}_*(f_m(x_m),
f_m(y_m))=\widetilde{h}_*(f_m(A)) \geqslant \delta.$ We may assume
that $x_m \rightarrow x_0,$ $y_m\rightarrow y_0$ as
$m\rightarrow\infty,$ $x_0, y_0\in A.$ If $\partial
D\ne\varnothing,$ we put $\varepsilon_0(x):=\min\{\widetilde{h}(A,
\partial D), r_*, r_0(x)\},$ where $r_*$ is defined in
Proposition~\ref{pr2} and $r_0(x)$ is a number from the conditions.
Otherwise, we set $\varepsilon_0(x)=\min\{r_*, r_0(x)\}.$ Let
$z_m\in A$ be a point such that
\begin{equation}\label{eq2A_1}
\widetilde{h}_*(f_m(A), \partial f_m(D))=\widetilde{h}_*(f_m(z_m),
\zeta_m)\rightarrow 0\,,\qquad m\rightarrow\infty\,,
\end{equation}
for some $\zeta_m\in
\partial f_m(D).$ We also may consider that $z_m\rightarrow z_0\in
A$ as $m\rightarrow\infty.$ Now we sequently join the points $x_0,$
$z_0$ and $y_0$ by a rectifiable path $\gamma$ in $D.$ We cover
$|\gamma|$ with the balls $\widetilde{B}(x, \varepsilon_0(x)/4),$
$x\in|\gamma|.$ By the Heine-Borel-Lebesgue lemma, we may choose a
finite subcover $\widetilde{B}(w_i, \varepsilon_i/4),$ $1\leqslant
i\leqslant M_0,$ $\varepsilon_i:=\varepsilon_0(w_i).$ We may
consider that $x_0\in \widetilde{B}(w_1, \varepsilon_1/4).$ Let
$\gamma:[0, 1]\rightarrow D,$ $\gamma(0)=x_0,$ $\gamma(1/2)=z_0$ and
$\gamma(1)=y_0.$ There are two possibilities: either
$|\gamma|\subset \widetilde{B}(w_1, \varepsilon_1/4),$ or
$|\gamma|\setminus \widetilde{B}(w_1,
\varepsilon_1/4)\ne\varnothing.$ In the first case, we may consider
that $\gamma$ is geodesic path (that is $\pi(\gamma_*),$ while
$\gamma_*$ is the corresponding hyperbolic geodesic path in the
hyperbolic disk $B_h(0, \varepsilon_1/4)$ and $\pi$ is the universal
covering of ${\Bbb S}$). Otherwise, if $|\gamma|\setminus
\widetilde{B}(w_1, \varepsilon_1/4)\ne\varnothing,$ then by
\cite[Theorem~1.I.5.46]{Ku$_2$}, there is $t_1\in (0, 1)$ such that
$\gamma(t_1)\cap \widetilde{S}(w_1, \varepsilon_1/4).$ We may
consider that $t_1=\sup\limits_{t\in [0, 1]:\gamma(t)\in
\overline{\widetilde{B}(w_1, \varepsilon_1/4)}}t.$ Now, we join
$x_0$ and $\gamma(t_1)$ with the geodesic path $\gamma_1.$ We may
consider that $\gamma(t_1)\in \widetilde{B}(w_2, \varepsilon_2/4).$
Arguing similarly to the above, either $\left|\gamma|_{[t_1,
1]}\right|\subset \widetilde{B}(w_2, \varepsilon_2/4),$ or
$\left|\gamma|_{[t_1, 1]}\right|\setminus \widetilde{B}(w_2,
\varepsilon_2/4)\ne\varnothing.$ In the first case, we may consider
that $\gamma|_{[t_1, 1]}$ is a geodesic path. Otherwise, if
$\left|\gamma|_{[t_1, 1]}\right|\setminus \widetilde{B}(w_2,
\varepsilon_2/4)\ne\varnothing,$ then by
\cite[Theorem~1.I.5.46]{Ku$_2$}, there is $t_2\in (t_1, 1)$ such
that $\gamma(t_2)\cap \widetilde{S}(w_2, \varepsilon_2/4).$ We may
consider that $t_2=\sup\limits_{t\in [t_1, 1]:\gamma(t)\in
\overline{\widetilde{B}(w_2, \varepsilon_2/4)}}t.$ Now, we join
$\gamma(t_1)$ and $\gamma(t_2)$ with the geodesic path $\gamma_2.$
We may consider that $\gamma(t_2)\in \widetilde{B}(w_3,
\varepsilon_3/4).$ Etc. By repeating this procedure a finite number
of times, we obtain a new path $\widetilde{\gamma},$ sequently
joining $x_0, z_0$ and $y_0$ that is a piecewise geodesic path.
Therefore, without loss of generality, we assume that this is
precisely the original path $\gamma.$

\medskip We also may assume that $x_m\in \widetilde{B}(x_0, \varepsilon_0(x_0)/4),$
$y_m\in \widetilde{B}(y_0, \varepsilon_0(y_0)/4)$ and $z_m\in
\widetilde{B}(z_0, \varepsilon_0(z_0)/4)$ for all $m\in {\Bbb N}.$
Due to the above, we may consider that the points $x_0,$ $z_0$ and
$y_0$ are joined by the piecewise geodesic path $\gamma,$ containing
some points $x_0=w_0, w_1, w_2,\ldots w_k, w_{k+1}=z_0, w_{k+2},
w_{k+3},\ldots , w_p=y_0$ such that $\widetilde{h}(w_{i},
w_{i+1})\leqslant\varepsilon_i/4$ for any $0\leqslant i\leqslant
p-1.$ Observe that, by the triangle inequality,
\begin{gather}\delta\leqslant \widetilde{h}_*(f_m(x_m), f_m(y_m))\leqslant
\widetilde{h}_*(f_m(x_m), f_m(x_0))+\sum\limits_{i=1}^{k+1}
\widetilde{h}_*(f_m(w_{i-1}), f_m(w_i))\nonumber \\
\label{eq5D} +2\widetilde{h}_*(f_m(z_0),
f_m(z_m))+\sum\limits_{i=k+2}^{p} \widetilde{h}_*(f_m(w_{i-1}),
f_m(w_i))+\widetilde{h}_*(f_m(y_0), f_m(y_m))\,.
\end{gather}
Since the sum~(\ref{eq5D}) has only $p+4$ terms, there is at least
one term greater than or equal to $\delta_*:=\frac{\delta}{p+4}>0$
for an infinite number of values of $m.$ In what follows, our
reasoning should be divided into two cases:
\begin{equation}\label{eq1_7}
\text{1)} \qquad \widetilde{h}_*(f_m(z_0), f_m(z_m))\geqslant
\delta_*>0
\end{equation}
for infinitely many $m\in{\Bbb N},$
\begin{equation}\label{eq1_8}
\text{2)} \qquad \widetilde{h}_*(f_m(z_0), f_m(z_m))\rightarrow
0\,,\qquad \text{as}\quad m\rightarrow\infty\,.
\end{equation}
Let us first consider {\bf case 1)} (the simpler one). So, assume
that $\widetilde{h}_*(f_m(z_0), f_m(z_m))\geqslant \delta_*>0$ for
infinitely many $m\in{\Bbb N}.$ Without loss of generality, we may
assume that the above holds for every $m\in {\Bbb N}.$ We now show
that
\begin{equation}\label{eq1_5}
\widetilde{h}_*(\partial(f_m(\widetilde{B}(z_0,
\varepsilon_0(z_0)/2))))\geqslant r_{**}>0 \end{equation}
for all $m=1,2,\ldots $ and some $r_{**}>0.$ We argue similarly to
the proof of~(\ref{eq1_2}). Otherwise, if~(\ref{eq1_5}) is not true,
we may consider that, there is $z_*\in \overline{D_*}$ such that,
for every $\sigma>0$ there is $m=m(\sigma)\in {\Bbb N}$ with
\begin{equation}\label{eq1_6}
f_m(\widetilde{S}(z_0, \varepsilon_0(z_0)/2))\subset
\widetilde{B}(z_*, \sigma)\,.
\end{equation}
Let $\sigma<\delta_*/2.$ Now, due to~(\ref{eq1_7}), it follows that
there is $\omega_m\in f_m(\widetilde{B}(z_0,
\varepsilon_0(z_0)/2))\setminus \widetilde{B}(z_*, \sigma).$ By the
assumption, $\partial f_m(D)\ne\varnothing.$ Let $\Delta_m\in
\partial f_m(D).$ Now, we join the points $\omega_m$ and $\Delta_m$
by a path $\beta_m$ in ${\Bbb S}_*\setminus \widetilde{B}(z_*,
\sigma).$ (The latter is possible due to the path connectedness of
$\widetilde{S}(z_*, \sigma)$ for sufficiently small $\sigma>0.$)
Since $|\beta_m|\cap f_m(\widetilde{B}(z_0,
\varepsilon_0(z_0)/2))\ne \varnothing \ne |\beta_m|\setminus
f_m(\widetilde{B}(z_0, \varepsilon_0(z_0)/2)),$ by
\cite[Theorem~1.I.5.46]{Ku$_2$} $|\beta_m|\cap
\partial f_m(\widetilde{B}(z_0,
\varepsilon_0(z_0)/2))\ne \varnothing.$ Since $\partial
f_m(\widetilde{B}(z_0, \varepsilon_0(z_0)/2))=f_m(\widetilde{S}(z_0,
\varepsilon_0(z_0)/2))$ and $\beta_m$ lies outside $
\widetilde{B}(z_*, \sigma),$ the latter contradicts
with~(\ref{eq1_6}). This contradiction proves~(\ref{eq1_5}). Now,
$f_m\in \frak{A}^{\delta_*, r_{**}}_{A_1, Q}(\widetilde{B}(z_0,
\varepsilon_0(z_0)/2))$ for $A_1:=\overline{\widetilde{B}(z_0,
\varepsilon_0(z_0)/4)}$ and by Lemma~\ref{lem2} we obtain that
$\widetilde{h}_*(f_m(A_1),
\partial f_m(\widetilde{B}(z_0, \varepsilon_0(z_0)/2)))$ $\geqslant
\delta_1>0$ for some $\delta_1>0.$ Since $z_m\in A_1$ for
sufficiently large $m\in {\Bbb N},$ it follows that
$\widetilde{h}_*(f_m(z_m),
\partial f_m(\widetilde{B}(z_0, \varepsilon_0(z_0)/2)))$ $\geqslant \delta_1>0$ and,
consequently, by Proposition~\ref{pr6} $\widetilde{h}_*(f_m(z_m),
\partial f_m(D))$ $\geqslant \delta_2>0$ for some $\delta_2>0.$
However, since $\widetilde{h}_*(f_m(A), \partial
f_m(D))=\widetilde{h}_*(f_m(z_m), \zeta_m)$ for some $\zeta_m\in
\partial f_m(D),$ we have that $\widetilde{h}_*(f_m(z_m), \zeta_m)\geqslant
\widetilde{h}_*(f_m(z_m), \partial f_m(D))\geqslant \delta_2>0$ for
sufficiently large $m\in {\Bbb N}.$ The latter contradicts
with~(\ref{eq1_9}).

\medskip
Let us consider the {\bf case 2),} when $\widetilde{h}_*(f_m(z_0),
f_m(z_m))\rightarrow 0$ as $m\rightarrow\infty.$ Here there are 2)
subcases:

\medskip
$2_1)$ $\widetilde{h}_*(f_m(z_0), f_m(z_m))\rightarrow 0$ as
$m\rightarrow\infty,$ $\widetilde{h}_*(f_m(w_{i-1}),
f_m(w_i))\rightarrow 0$ as $m\rightarrow\infty$ for any $0\leqslant
i\leqslant p,$ however, one of two following conditions hold: either
$\widetilde{h}_*(f_m(x_m), f_m(x_0))\geqslant \delta_*>0$ or
$\widetilde{h}_*(f_m(y_m), f_m(y_0))\geqslant \delta_*>0$ for
infinitely many $m\in {\Bbb N};$

\medskip
$2_2)$ $\widetilde{h}_*(f_m(z_0), f_m(z_m))\rightarrow 0$ as
$m\rightarrow\infty,$ besides that, $\widetilde{h}_*(f_m(w_{i-1}),
f_m(w_i))\geqslant \delta_*>0$ for infinitely many $m\in {\Bbb N}$
and at least for one $0\leqslant i\leqslant p.$

\medskip
Let us consider {\bf the situation $2_1)$} firstly. Assume that
$\widetilde{h}_*(f_m(x_m), f_m(x_0))\geqslant \delta_*>0$ for
infinitely many $m\in {\Bbb N}$ (another case, when
$\widetilde{h}_*(f_m(y_m), f_m(y_0))\geqslant \delta_*>0,$ may be
consider similarly). Now, we have that
$\widetilde{h}_*(f_{m_k}(x_{m_k}), f_{m_k}(x_0))\geqslant
\delta_*>0$ for all $k\in {\Bbb N}.$ We may consider that the latter
holds for $m$ instead $m_k.$

Arguing similarly to the point~1), we may obtain that $$
\widetilde{h}_*(\partial(f_m(\widetilde{B}(x_0,
\varepsilon_0(x_0)/2))))\geqslant r_{**}>0$$
for some $r_{**}>0,$ while $\widetilde{h}_*(f_m(x_0),
f_m(x_m))\geqslant \delta_*>0$ for $m\in{\Bbb N}$ and, consequently,
$\widetilde{h}_*(f_m(A_1))\geqslant \delta_*>0,$ where
$A_1:=\overline{\widetilde{B}(x_0, \varepsilon_0(x_0)/4)}.$ Now,
$f_m\in \frak{A}^{\delta_*, r_{**}}_{A_1, Q}(\widetilde{B}(x_0,
\varepsilon_0(x_0)/2))$ and by Lemma~\ref{lem2} we obtain that
$$\widetilde{h}_*(f_m(A_1), \partial f_m(\widetilde{B}(x_0,
\varepsilon_0(x_0)/2)))\geqslant \delta_1>0$$ for some $\delta_1>0.$
Since $x_m\in A_1,$ it follows that $\widetilde{h}_*(f_m(x_m),
\partial f_m(\widetilde{B}(x_0, \varepsilon_0(x_0)/2)))\geqslant \delta_1>0$ and,
consequently, by Proposition~\ref{pr6} $\widetilde{h}_*(f_m(x_m),
\partial f_m(D))\geqslant \delta_2>0$ for some $\delta_2>0.$

Observe that $\widetilde{h}_*(f_m(x_m),
\partial f_m(D))\leqslant \widetilde{h}_*(f_m(x_m), \zeta_m),$ where $\zeta_m\in
\partial f_m(D)$ is from~(\ref{eq2A_1}). Now
\begin{equation}\label{eq3_1}
\widetilde{h}_*(f_m(x_m), \zeta_m)\geqslant \delta_2>0\,,
m=1,2,\ldots \,,
\end{equation}
Since $\widetilde{h}_*(f_m(w_{i-1}), f_m(w_i))\rightarrow 0$ as
$m\rightarrow\infty$ for any $0\leqslant i\leqslant p$ and
$\widetilde{h}_*(f_m(z_0), f_m(z_m))\rightarrow 0$ as
$m\rightarrow\infty,$ by the triangle inequality and~(\ref{eq3_1})
we obtain that
\begin{gather}\widetilde{h}_*(f_m(z_m), \zeta_m)\geqslant
\widetilde{h}_*(\zeta_m, f_m(x_m))-\widetilde{h}_*(f_m(x_m),
f_m(z_m))\nonumber\\
\label{eq4_1} \geqslant \widetilde{h}_*(\zeta_m,
f_m(x_m))-\sum\limits_{i=1}^{k+1} \widetilde{h}_*(f_m(w_{i-1}),
f_m(w_i))-\widetilde{h}_*(f_m(z_0), f_m(z_m))\geqslant \delta_2/2>0
\end{gather}
since the latter two terms with a minus sign tend to zero by the
assumption in~$2_1).$ The relation~(\ref{eq4_1}) contradicts
with~(\ref{eq2A_1}). The latter completes the consideration of
item~$2_1).$

\medskip
It remains to consider the situation~$2_2).$ Let
$\widetilde{h}_*(f_m(w_{i-1}), f_m(w_i))\geqslant \delta_*>0$ for
infinitely many $m\in {\Bbb N}$ and at least for one $0\leqslant
i\leqslant p.$ We may consider that the latter holds for all $m\in
{\Bbb N}.$ Now we consider the case $0\leqslant i\leqslant k+1,$
while the consideration of the case $k+1\leqslant i\leqslant p$ may
be done similarly. We also may consider that $i$ be the largest
index between $0\leqslant i\leqslant p$ for which the relation
$\widetilde{h}_*(f_m(w_{i-1}), f_m(w_i))\geqslant \delta_*>0$ holds.
Now, $\widetilde{h}_*(f_m(w_{i-1}), f_m(w_i))\geqslant \delta_*>0$
for $m\in {\Bbb N},$ and $\widetilde{h}_*(f_m(w_{l-1}),
f_m(w_l))\rightarrow 0$ as $m\rightarrow\infty$ for all
$i<l\leqslant k+1.$

\medskip
We put $A_1:=\overline{B(w_{i}, \varepsilon_i/4)}.$ Now,
$\widetilde{h}_*(f_m(A_1))\geqslant\delta_*>0$
for all $m=1,2,\ldots ,$ because $w_{i-1}, w_{i}\in A_1$ for every
$m\in {\Bbb N}.$ Arguing similarly to the proof of~(\ref{eq1_5}), we
may prove that $\widetilde{h}_*(\partial(f_m(\widetilde{B}(w_{i},
\varepsilon_i/2))))\geqslant r_{**}>0$ for all $m=1,2,\ldots $ and
some $r_{**}>0$ while $\partial(f_m(\widetilde{B}(w_{i},
\varepsilon_i/2)))$ is connected as a image of the connected set
$\widetilde{S}(w_{i}, \varepsilon_i/2)$ under the mapping $f_m.$
Now, $f_m\in \frak{A}^{\delta_*, r_{**}}_{A_1,
Q}(\widetilde{B}(w_{i}, \varepsilon_i/2))$ and by Lemma~\ref{lem2}
$\widetilde{h}_*(f_m(A),
\partial f_m(\widetilde{B}(w_{i}, \varepsilon_i/2)))\geqslant \delta_1>0$ for
some $\delta_1>0.$ Since $w_{i}\in A_1,$ it follows that
$$\widetilde{h}_*(f_m(w_{i}),
\partial f_m(\widetilde{B}(w_{i}, \varepsilon_i/2)))\geqslant \delta_1>0$$ and,
consequently, by Proposition~\ref{pr6} we have that
$\widetilde{h}_*(f_m(w_{i}),
\partial f_m(D))\geqslant \delta_2>0$ for some $\delta_2>0.$
Observe that $$\widetilde{h}_*(f_m(w_{i}),
\partial f_m(D))\leqslant \widetilde{h}_*(f_m(w_{i}), \zeta_m)\,,$$ where $\zeta_m\in
\partial f_m(D)$ is from~(\ref{eq2A_1}). Now
\begin{equation}\label{eq3B}
\widetilde{h}_*(f_m(w_{i}), \zeta_m)\geqslant \delta_2>0\,,
m=1,2,\ldots .
\end{equation}
Since $\widetilde{h}_*(f_m(w_{l-1}), f_m(w_l))\rightarrow 0$ as
$m\rightarrow\infty$ for any $i<l\leqslant k+1$ and simultaneously
$\widetilde{h}_*(f_m(z_0), f_m(z_m))\rightarrow 0$ as
$m\rightarrow\infty,$ by the triangle inequality and~(\ref{eq3B}) we
obtain that
\begin{gather}\widetilde{h}_*(f_m(z_m), \zeta_m)\geqslant
\widetilde{h}_*(\zeta_m, f_m(w_{i}))-\widetilde{h}_*(f_m(w_{i}),
f_m(z_m))\nonumber \\
\label{eq4A} \geqslant \widetilde{h}_*(\zeta_m,
f_m(w_{i}))-\sum\limits_{l=i+1}^{k+1} \widetilde{h}_*(f_m(w_{l-1}),
f_m(w_l))-\widetilde{h}_*(f_m(z_0), f_m(z_m))\geqslant \delta_2/2>0
\end{gather}
because the latter two terms with a minus sign tend to zero by the
assumption in this subparagraph~$2_2).$ The relation~(\ref{eq4A})
contradicts with~(\ref{eq2A_1}). The latter completes the
consideration of item~$2_2).$ The proof is complete.
\end{proof}

\begin{corollary}\label{cor5_1}
{\, The statement of Lemma~\ref{lem4} remains true if, instead of
the condition regarding the integrability of the function $Q$ over
spheres with respect to some set $E_1$ is replaced by a simpler
condition: $Q\in L_{\rm loc}^1(D).$}
\end{corollary}

{\it Proof of Theorem~\ref{th2}.} Let $K$ be a compactum in $D.$
Observe that, $\partial f(D)\ne \varnothing$ for
$f\in\frak{B}^{\delta}_{A, Q}(D, D_*),$ because $\partial
D\ne\varnothing$ and $f$ is a homeomorphism (see, e.g.,
\cite[Proposition~13.5]{MRSY}). By Lemma~\ref{lem4} and
Proposition~\ref{pr4}, $\frak{B}^{\delta}_{A, Q}(D, D_*)\subset
\frak{F}^{\delta_1}_{K, Q}(D, D_*)$ for some $\delta_1>0.$ The
desired conclusion follows now by Theorems~\ref{th1}
and~\ref{th2}.~$\Box$

\medskip
{\it Proof of Corollaries~\ref{cor6} and~\ref{cor7}} is similar to
the proofs of of Corollaries~\ref{cor1} and~\ref{cor2} and,
therefore me be omitted.~$\Box$

\section{The convergence theorem}

Similarly to \cite{RS$_1$}, we have the following, cf.
\cite[Theorem~1]{Cr$_1$}, \cite[Theorem~1.2]{Cr$_2$} and
\cite[Theorem~4.2]{RS$_2$}.

\begin{theorem}\label{th4}
{\, Let $D$ and $D_*$ be domains in ${\Bbb S}$ and ${\Bbb S}_*,$
respectively, besides that, $\overline{D_*}$ is a compactum in
${\Bbb S}_*.$ Assume that, for each point $p_0\in D$ there is
$r_0=r_0(p_0)>0,$ $r_0<\sup\limits_{p\in D}\widetilde{h}(p, p_0),$
such that for any $0<r_1<r_2<r_0$ there is a set $E_1\subset[r_1,
r_2]$ of a positive linear Lebesgue measure such that the function
$Q$ is integrable over $\widetilde{S}(p_0, r)$ for every $r\in E_1.$
Let $f_m:D\rightarrow D_*,$ $m=1,2,\ldots, $ be a sequence of
homeomorphisms satisfying the relations
(\ref{eq3*!!A})--(\ref{eq8BC}) for all $p_0\in D$ and all
$0<r_1<r_2<r_0$ and some $r_0>0,$ $r_0<\widetilde{h}(p_0, \partial
D),$ while $\widetilde{h}(p_0, \partial D):=\infty$ if $\partial
D=\varnothing.$ If $f_m$ converges to $f:D\rightarrow D_*$ locally
uniformly, then $f$ either is a homeomorphism $f:D\rightarrow D_*$
or a constant $c\in \overline{D_*}.$}
\end{theorem}

\begin{proof}
If $f\equiv const,$ there is nothing to prove, so that we assume
that $f\not\equiv const.$ We firstly prove that $f$ is discrete.
Assume the contrary. Then there is $x_0 \in D$ and a sequence
$x_m\in D,$ $m=1,2,\ldots ,$ $x_m \ne x_0,$ such that
$x_m\rightarrow x_0$ as $m\rightarrow \infty$ and $f(x_m)=f(x_0).$
Replacing $D$ by some neighborhood $U$ of $x_0$ and reducing this
neighborhood, if required, we my consider that $\partial D$ contains
infinitely many points. Observe that $E_0=\{x\in D: f(x)=f(x_0)\}$
is closed in $D$ by the continuity of $f$ and does not coincide with
$D,$ because $f\not \equiv const.$ Thus, we may assume that $x_0$
may be replaced by a non isolated boundary point of $E_0.$

\medskip
Now, there are $\varepsilon_0>0,$ $0<\varepsilon_0<
\frac{1}{2}\cdot\widetilde{h}(x_0,
\partial D),$ and $a, b\in \widetilde{B}(x_0, \varepsilon_0),$
$a\ne b,$ such that $\widetilde{h}_*(f(a), f(b)):=\delta>0.$ Since
$f_m$ converges to $f$ locally uniformly, $\widetilde{h}_*(f_m(a),
f_m(b))\geqslant \delta/2>0$ for sufficiently large $m\in {\Bbb N}.$
Let $A:=|\gamma|,$ where $\gamma$ is a path joining $a$ and $b$ in
$\widetilde{B}(x_0, \varepsilon_0).$ Now,
$\widetilde{h}_*(f_m(A))\geqslant \delta>0$ for all $m=1,2,\ldots.$

\medskip
Thus, $f_m\in \frak{B}^{\delta}_{A, Q}(D, D_*).$ Let $K$ be some
closed ball centered at $x_0$ in entirely lying in $D.$ By
Theorem~\ref{th2} the family $\frak{B}^{\delta}_{A, Q}(D, D_*)$ is
uniformly light on $K,$ that is, given $\varepsilon>0$ there is
$\delta=\delta(\varepsilon, K)>0$ such that $\widetilde{h}_*(f_m(x),
f_m(x_0))\geqslant \delta$ whenever $\widetilde{h}(x, x_0)\geqslant
\varepsilon>0.$ Taking here the limit as $m\rightarrow\infty,$ we
obtain that $\widetilde{h}_*(f(x), f(x_0))\geqslant \delta.$ This
means that $f$ is discrete in $D.$

\medskip
We now prove that $f$ is open. Let $A$ be an open set and let
$x_0\in A.$ Since $A$ is open, there is $\varepsilon_0>0$ such that
$\overline{B(x_0, \varepsilon_0)}\subset A$ for some
$\varepsilon_0>0.$ We may assume that $\overline{B(x_0,
\varepsilon_0)}$ is a compactum in ${\Bbb S}.$ By
Corollary~\ref{cor7} the family $f_m$ is uniformly open in
$\overline{\widetilde{B}(x_0, \varepsilon_0)},$ i.e., there is
$r_0>0$ such that
\begin{equation}\label{eq1AA}
\widetilde{B}(f_m(x_0), r_0)\subset f_m(\widetilde{B}(x_0,
\varepsilon_1))
\end{equation}
for every $m\in {\Bbb N}$ and any $0<\varepsilon_1<\varepsilon_0.$
Let $y\in \widetilde{B}(f(x_0), r_0/2).$ By the uniform convergence
of $f_m$ to $f$ and by the triangle inequality we obtain that
$\widetilde{B}(f(x_0), r_0/2)\subset \widetilde{B}(f_m(x_0), r_0).$
Now, by~(\ref{eq1AA}), $y\in f_m(\widetilde{B}(x_0,
\varepsilon_1)).$ Thus, $y=f_m(x_m)$ for some $x_m\in
\widetilde{B}(x_0, \varepsilon_1).$ Due to the compactness of
$\overline{\widetilde{B}(x_0, \varepsilon_1)},$ we may consider that
$x_m\rightarrow z_0\in \overline{\widetilde{B}(x_0, \varepsilon_1)}$
as $m\rightarrow\infty.$ By the continuity of $f$ in $A,$ since
$\overline{\widetilde{B}(x_0, \varepsilon_1)}\subset A,$ we obtain
that $f(x_m)\rightarrow f(z_0)$ as $m\rightarrow\infty.$ So, we have
that $f(x_m)\rightarrow f(z_0)$ as $m\rightarrow\infty$ and
simultaneously $y=f_m(x_m)$ for sufficiently large $m\in {\Bbb N}.$
Thus, by the triangle inequality,
\begin{gather*}\widetilde{h}_*(y, f(z_0))=\widetilde{h}_*(f_m(x_m),
f(z_0))\\
\leqslant \widetilde{h}_* (f_m(x_m), f(x_m))+
\widetilde{h}_*(f(x_m), f(z_0))\rightarrow 0\,, \end{gather*}
$m\rightarrow\infty.$ Thus, $y=f(z_0)\in
f(\overline{\widetilde{B}(x_0, \varepsilon_1)})\subset f(A).$ So,
$y\in f(A),$ i.e., $\widetilde{B}(f(x_0), r_0/2)\subset f(A).$ The
latter means that $f(A)$ is open, as required.

\medskip
Finally, we show that the mapping $f$ is injective, cf.
\cite[Theorem~4]{IRS}. We argue by the contradiction, i.e., assume
that there are $x_1, x_2\in D$ such that $x_1\ne x_2,$ however,
$f(x_1)=f(x_2).$ Set $B_t:=\widetilde{B}(x_1, t).$ Due to the
discreteness of $f,$ there is $t_0>0$ such that $x_2\not\in
\overline{\widetilde{B}(x_1, t)}$ for any $t\in [0, t_0].$

\medskip
Let $y_m:=f_m(x_1),$ $z_m=f_m(x_2)$ and $C_m:=f_m(B_t).$ Since $f_m$
converges to $f$ locally uniformly to $f,$ we have that
$y_m\rightarrow f(x_1)$ and $z_m\rightarrow f(x_2)$ as
$m\rightarrow\infty,$ so that $y_m$ lies in some normal neighborhood
$U_0$ of $f(x_1)$ for sufficiently large $m\in {\Bbb N}.$ Let
$\sigma>0$ be such that $\widetilde{B}(f(x_1), \sigma)\subset U_0.$
Observe that $\partial C_m=f_m(\partial B_t)$ are compacta as the
continuous images of compacta. Therefore, there is $x_{m,t} \in
\partial B_t$ such that
\begin{equation}\label{eq3B_1}
\widetilde{h}_*(y_m, \partial C_m) = \widetilde{h}_*(y_m,
f_m(x_{m,t}))=\widetilde{h}_*(f_m(x_1), f_m(x_{m,t}))\,,\qquad
m=1,2,\ldots\,.
\end{equation}
Given $\varepsilon>0.$ Now, by~(\ref{eq3B_1}) and by the triangle
inequality, we obtain that
$$\widetilde{h}_*(y_m, \partial C_m)\leqslant
\widetilde{h}_*(f_m(x_1), f(x_1))+\widetilde{h}_*(f(x_1),
f(x_{m,t}))+\widetilde{h}_*(f(x_{m,t}), f_m(x_{m,t}))\,.$$
Since $f_m$ converges to $f$ locally uniformly, the latter
inequality implies that there are
$t^{\,*}_0=t^{\,*}_0(\varepsilon)>0$ and $m_0=m_0(\varepsilon)\in
{\Bbb N}$ such that
\begin{equation}\label{eq12A}
\widetilde{h}_*(y_m, \partial C_m)<\varepsilon\,,\qquad m\geqslant
m_0\,,\qquad t\in [0,t^{\,*}_0]\,.
\end{equation}
We may consider that $0<t_0<t^{\,*}_0.$ Let $y\in \widetilde{B}(y_m,
\widetilde{h}_*(y_m,
\partial C_m)).$ Now, by the triangle inequality and~(\ref{eq12A})
\begin{gather*}\widetilde{h}_*(y, f(x_1))\leqslant \widetilde{h}_*(y, y_m)+
\widetilde{h}_*(y_m, f(x_1)) \leqslant \widetilde{h}_*(y_m,
\partial C_m)+\widetilde{h}_*(y_m, f(x_1))<2\varepsilon
\end{gather*}
for all $m\geqslant m_1,$ $t\in [0, t^{\,*}_0]$ and some $m_1>m_0.$
If we take $2\varepsilon<\sigma,$ we have from here that
$$\widetilde{B}(y_m, \widetilde{h}_*(y_m,
\partial C_m))\subset \widetilde{B}(f(x_1),
\sigma)\subset U_0\,.$$ The latter means that $\widetilde{B}(y_m,
\widetilde{h}_*(y_m,
\partial C_m))\subset U_0$ for $0<t<t_0$ and $m\geqslant m_1.$ Now, $\widetilde{B}(y_m,
\widetilde{h}_*(y_m,
\partial C_m))$ is connected because it belongs to the normal
neighborhood of some point for sufficiently large $m$ and all
$0<t<t_0$ and, therefore, $\widetilde{B}(y_m, \widetilde{h}_*(y_m,
\partial C_m))$ is isometrical to some hyperbolic disk in ${\Bbb D}.$
We now show that
\begin{equation}\label{eq1_*}
\widetilde{B}(y_m, \widetilde{h}_*(y_m, \partial C_m))\subset
C_m\,,\qquad m=1,2,\ldots\,.
\end{equation}
Assume the contrary. Then is $y\in \widetilde{B}(y_m,
\widetilde{h}_*(y_m,
\partial C_m))\setminus C_m.$ Now, $$\widetilde{B}(y_m, \widetilde{h}_*(y_m,
\partial C_m))\setminus C_m\ne\varnothing\ne C_m\cap \widetilde{B}(y_m,
\widetilde{h}_*(y_m,
\partial C_m))$$ because $y_m\in C_m,$ as well.
Since by the proving above, $\widetilde{B}(y_m, \widetilde{h}_*(y_m,
\partial C_m))$ is connected for sufficiently large $m=1,2,\ldots $,
by \cite[Theorem~1.I.5.46]{Ku$_2$} we obtain that
$\widetilde{B}(y_m, \widetilde{h}_*(y_m,
\partial C_m))\cap \partial  C_m\ne\varnothing.$ Now, there is $z\in \widetilde{B}(y_m,
\widetilde{h}_*(y_m,
\partial C_m))\cap \partial C_m.$ Thus,
$$\widetilde{h}_*(y_m,
\partial C_m)\leqslant \widetilde{h}_*(z, y_m)<\widetilde{h}_*(y_m,
\partial C_m)\,.$$
The contradiction obtained above proves~(\ref{eq1_*}). Now,
by~(\ref{eq1_*})
\begin{equation}\label{eq2C}
\overline{\widetilde{B}(y_m, \widetilde{h}_*(y_m, \partial
C_m))}\subset \overline{C_m}\,,\qquad m=1,2,\ldots\,.
\end{equation}
We show that
\begin{equation}\label{eq5_1}
\widetilde{h}_*(y_m, \partial C_m)<\widetilde{h}_*(y_m,
z_m)\,,\qquad m=1,2,\ldots\, .
\end{equation}
In the contrary case we obtain that $\widetilde{h}_*(y_m,
z_m)\leqslant \widetilde{h}_*(y_m,
\partial C_m)$ for some $m\in{\Bbb N},$ consequently we obtain that
$z_m\in \overline{\widetilde{B}(y_m, \widetilde{h}_*(y_m, \partial
C_m))}.$ But now, by~(\ref{eq2C}) $z_m\in \overline{C_m}$ that
contradicts the choice of $z_m=f_m(x_2)$ and that $f_m$ is a
homeomorphism.

\medskip

By~(\ref{eq5_1}) and~(\ref{eq3B_1}),
$$\widetilde{h}_*(y_m,
f_m(x_{m,t}))\leqslant \widetilde{h}_*(y_m, z_m)\,,\qquad
m=1,2,\ldots\,,$$
or, equivalently,
\begin{equation}\label{eq6_1}
\widetilde{h}_*(f_m(x_1), f_m(x_{m,t}))\leqslant
\widetilde{h}_*(f_m(x_1), f_m(x_2))\,,\qquad m=1,2,\ldots\,.
\end{equation}
By the compactness of $\partial B_t$, for every $t\in (0, t_0],$
there is $x_t\in \partial B_t$ such that $\widetilde{h}_*(x_{m_k,t},
x_t)\rightarrow 0$ as $k\rightarrow \infty$ for some subsequence
$m_k,$ $k=1,2,\ldots .$ Since the locally uniform convergence of
continuous functions in a metric space implies the continuous
convergence (see Theorem 21.X.3 in \cite{Ku$_1$}), we have that
$\widetilde{h}_*(f_{m_k}(x_{m_k,t}),f(x_{t})) \rightarrow 0$ as
$k\rightarrow \infty.$ Consequently, by~(\ref{eq6_1}) we obtain that
$$\widetilde{h}_*(f(x_1),f(x_t))\leqslant \widetilde{h}_*(f(x_1),f(x_2))$$
for any $t\in(0, t_0].$
However, since $f(x_1)=f(x_2)$ we obtain that $f(x_t)=f(x_1)$ for
every $t\in (0,t_0].$ The latter contradicts to the discreteness of
$f.$ Thus, $f$ is injective. In addition, $f^{\,-1}$ is continuous
in $f(D)$ because $f$ is open by the proving above.
\end{proof}

\section{Consequences for Sobolev classes}

{\it A (maximal) dilatation} of $f$ at $z$ is defined in local
coordinates by the relation
\begin{equation*}
K_f(z)=\frac{|f_z|+|f_{\overline{z}}|}{|f_z|-|f_{\overline{z}}|}
\end{equation*}
for $J_f(z)\ne 0,$ $K_f(z)=1$ for $\Vert f^{\,\prime}(z)\Vert=0$ and
$K_f(z)=\infty$ otherwise. It is not difficult to see that, $K_f$
does not depend on local coordinates because the transition mappings
between two charts are conformal by the definition of the Riemannian
surface.

\medskip
Let $D$ and $D_*$ be domains in ${\Bbb S}$ and ${\Bbb S}_*,$
correspondingly. A mapping $f:D\rightarrow D_{\,*}$ is called a {\it
mapping with finite distortion,} if $f\in W_{\rm loc}^{1, 1}(D)$
and, in addition, there is a Lebesgue measurable function $K(z),$
$K(z)<\infty$ a.e., such that $\Vert f^{\,\prime}(z)\Vert^2\leqslant
K(z)\cdot J_f(z)$ for almost all $z \in D.$ (Here $\Vert
f^{\,\prime}(z)\Vert$ and $J_f(z)$ denote the matrix norm of
$f^{\,\prime}$ and the jacobian of $f,$ respectively; see, e.g.,
\cite{RV$_1$}). Due to the Lemma~3.1 in \cite{RV$_1$}, the following
statement is true.

\begin{proposition}\label{pr5}
{\, Let $D$ and $D_*$ be domains on Riemannian surfaces ${\Bbb S}$
and ${\Bbb S}_*.$ If $f:D\rightarrow D_*$ is a homeomorphism of
finite distortion with $K_f\in L^1_{\rm loc},$ then $f$ satisfies
(\ref{eq3*!!A})--(\ref{eq8BC}) for any $p_0\in \overline{D}$ with
$Q=K_f.$}\end{proposition}

Given domains $D$ and $D_*$ in ${\Bbb S}$ and ${\Bbb S}_*,$
respectively, a function $Q:D\rightarrow [0, \infty]$ measurable by
the measure $\widetilde{v},$ a compact set $K\subset D$ and a number
$\delta>0$ denote by $\frak{Sob^1}^{\delta}_{K, Q}(D, D_*)$ a family
of homeomorphisms $f:D\rightarrow D_*$ such that $K_f\in L^1_{\rm
loc}$ and $K_f(x)\leqslant Q(x)$ a.e., and, besides that,
$\widetilde{h}_*(f(K),
\partial f(D))=\inf\limits_{x\in f(K), y\in \partial f(D)}
\widetilde{h}_*(x, y)\geqslant \delta.$

\medskip
Given domains $D, D_*$ in ${\Bbb S}$ and ${\Bbb S}_*,$ respectively,
a non-degenerate continuum $A\subset D,$ a measurable function
$Q:D\rightarrow [0, \infty]$ and a number $\delta>0$ we denote by
$\frak{Sob^2}^{\delta}_{A, Q}(D, D_*)$ a family of all
homeomorphisms $f:D\rightarrow D_*$ such that $K_f\in L^1_{\rm loc}$
and $K_f(x)\leqslant Q(x)$ a.e., and, besides that,
$\widetilde{h}_*(f(A))\geqslant \delta.$

\medskip
The following statements directly follows by Theorems~\ref{th1},
\ref{th2} and Corollaries~\ref{cor1}, \ref{cor2}, \ref{cor6} and
\ref{cor7} taking into account Proposition~\ref{pr5}.

\begin{theorem}\label{th1_1}
{\,Let $D$ and $D_*$ be domains in ${\Bbb S}$ and ${\Bbb S}_*,$
respectively, while $D$ has at least two different boundary points
$a$ and $b$ and, besides that, $\overline{D_*}$ is a compactum in
${\Bbb S}_*.$ Assume that, for each point $p_0\in D$ there is
$r_0=r_0(p_0)>0,$ $r_0<\sup\limits_{p\in D}\widetilde{h}(p, p_0),$
such that for any $0<r_1<r_2<r_0$ there is a set $E_1\subset[r_1,
r_2]$ of a positive linear Lebesgue measure such that the function
$Q$ is integrable over $\widetilde{S}(p_0, r)$ for every $r\in E_1.$
Then the family $\frak{Sob^1}^{\delta}_{K, Q}(D, D_*)$ is uniformly
light on $K;$ moreover, $\frak{Sob^1}^{\delta}_{K, Q}(D, D_*)$ is
equi-uniform on $K.$

\medskip
Moreover, if (\ref{eq3*!!A})--(\ref{eq8BC}) hold for any $p_0\in D,$
any $0<r_1<r_2<\sup\limits_{p\in D}\widetilde{h}(p, p_0),$ any $f\in
\frak{Sob^1}^{\delta}_{K, Q}(D, D_*)$ and $Q\in L^1(D),$ then there
exist $C_1>0$ such that
$$
\widetilde{h}_*(f(x), f(y))\geqslant
C_1\cdot\exp\left\{-\frac{2\Vert Q\Vert_1}{C{(\widetilde{h}(x,
y))}^2}\right\}
$$
for all $x, y\in K,$ $x\ne y,$ and every $f\in
\frak{Sob^1}^{\delta}_{K, Q}(D, D_*),$ where $C$ is a constant from
(\ref{eq3C}) for $n=2$ and $X={\Bbb C}.$ }
\end{theorem}

\begin{corollary}\label{cor1_1}
{\, The statement of Theorem~\ref{th1_1} remains true if we replace
the assumption on $Q$ mentioned above concerning the integrability
of $Q$ over circles $\widetilde{S}(p_0, r)$ by a simpler condition:
$Q\in L_{\rm loc}^1(D).$}
\end{corollary}

\begin{corollary}\label{cor2_1}{\bf \,(An analogue of Koebe-Bloch
theorem)}.~{ Under assumptions of Theorem~\ref{th1_1}, the family
$\frak{Sob^1}^{\delta}_{K, Q}(D, D_*)$ is uniformly open on $K.$ In
particular, the latter is true if we replace the assumption on $Q$
mentioned in Theorem~\ref{th1_1} concerning the integrability of $Q$
over circles $\widetilde{S}(p_0, r)$ by a simpler condition: $Q\in
L_{\rm loc}^1(D).$}
\end{corollary}

\begin{theorem}\label{th2_1}
{\,Let $D$ and $D_*$ be domains in ${\Bbb S}$ and ${\Bbb S}_*,$
respectively, while $D$ has at least two different boundary points
$a$ and $b$ and, besides that, $\overline{D_*}$ is a compactum in
${\Bbb S}_*.$ Assume that, for each point $p_0\in D$ there is
$r_0=r_0(p_0)>0,$ $r_0<\sup\limits_{p\in D}\widetilde{h}(p, p_0),$
such that for any $0<r_1<r_2<r_0$ there is a set $E_1\subset[r_1,
r_2]$ of a positive linear Lebesgue measure such that the function
$Q$ is integrable over $\widetilde{S}(p_0, r)$ for every $r\in E_1.$
Then the family $\frak{Sob^2}^{\delta}_{A, Q}(D, D_*)$ is uniformly
light on every compactum $K\subset D$ moreover,
$\frak{Sob^2}^{\delta}_{A, Q}(D, D_*)$ is equi-uniform on $K.$

\medskip
Moreover, if (\ref{eq3*!!A})--(\ref{eq8BC}) hold for any $p_0\in D,$
any $0<r_1<r_2<\sup\limits_{p\in D}\widetilde{h}(p, p_0),$ any $f\in
\frak{Sob^2}^{\delta}_{A, Q}(D, D_*)$ and $Q\in L^1(D),$ then there
exist $C_1>0$ such that the relation~(\ref{eq1D}) holds
for all $x, y\in K,$ $x\ne y,$ and every $f\in
\frak{Sob^2}^{\delta}_{A, Q}(D, D_*),$ where $C$ is a constant from
(\ref{eq3C}) for $n=2$ and $X={\Bbb C}.$ }
\end{theorem}

\begin{corollary}\label{cor6_1}
{\, The statement of Theorem~\ref{th2_1} remains true if we replace
the assumption on $Q$ mentioned above concerning the integrability
of $Q$ over circles $\widetilde{S}(p_0, r)$ by a simpler condition:
$Q\in L_{\rm loc}^1(D).$}
\end{corollary}

\begin{corollary}\label{cor7_1}{\bf \,(An analogue of Koebe-Bloch
theorem)}.~{ Under assumptions of Theorem~\ref{th2_1}, the family
$\frak{Sob^2}^{\delta}_{A, Q}(D, D_*)$ is uniformly open on every
compactum $K\subset D.$ In particular, the latter is true if we
replace the assumption on $Q$ mentioned in Theorem~\ref{th2_1}
concerning the integrability of $Q$ over circles $\widetilde{S}(p_0,
r)$ by a simpler condition: $Q\in L_{\rm loc}^1(D).$}
\end{corollary}

\section{On inverse Sobolev mappings}

Observe that, if $f:D\rightarrow {\Bbb R}^2,$ $D\subset {\Bbb R}^2,$
is a homeomorphism with a finite distortion, then $f^{\,-1}$ is a
homeomorphism with a finite distortion, as well (see
\cite[Theorem~1.3]{HK}). Now, we formulate the corresponding
statements for the inverse Sobolev classes which are direct
consequences of the results of the previous section.

\medskip
Given domains $D$ and $D_*$ in ${\Bbb S}$ and ${\Bbb S}_*,$
respectively, a function $Q:D\rightarrow [0, \infty]$ measurable by
the measure $\widetilde{v},$ a compact set $K\subset D$ and a number
$\delta>0$ denote by $\frak{Sob}^{-1\delta}_{K, Q}(D, D_*)$ a family
of homeomorphisms $g:D_*\rightarrow D,$ $g(D_*)=D,$ such that
$K_f\in L^1_{\rm loc}$ and $K_f(x)\leqslant Q(x)$ a.e.,
$f:=g^{\,-1},$ and $\widetilde{h}_*(f(K),
\partial f(D))=\inf\limits_{x\in f(K), y\in \partial f(D)}
\widetilde{h}_*(x, y)\geqslant \delta.$

\medskip
Given domains $D, D_*$ in ${\Bbb S}$ and ${\Bbb S}_*,$ respectively,
a non-degenerate continuum $A\subset D,$ a measurable function
$Q:D\rightarrow [0, \infty]$ and a number $\delta>0$ we denote by
$\frak{Sob}^{-2\delta}_{A, Q}(D, D_*)$ a family of all
homeomorphisms $g:D_*\rightarrow D,$ $g(D_*)=D,$ such that $K_f\in
L^1_{\rm loc}$ and $K_f(x)\leqslant Q(x)$ a.e., $f:=g^{\,-1},$ and
$\widetilde{h}_*(f(A))\geqslant \delta.$

\medskip
The following statements directly follows by Theorems~\ref{th1_1},
\ref{th2_1} and Corollaries~\ref{cor1_1}, \ref{cor2_1}, \ref{cor6_1}
and \ref{cor7_1} due to \cite[Theorem~1.3]{HK}.

\begin{theorem}\label{th1_2}
{\,Let $D$ and $D_*$ be domains in ${\Bbb S}$ and ${\Bbb S}_*,$
respectively, while $D$ has at least two different boundary points
$a$ and $b$ and, besides that, $\overline{D_*}$ is a compactum in
${\Bbb S}_*.$ Assume that, for each point $p_0\in D$ there is
$r_0=r_0(p_0)>0,$ $r_0<\sup\limits_{p\in D}\widetilde{h}(p, p_0),$
such that for any $0<r_1<r_2<r_0$ there is a set $E_1\subset[r_1,
r_2]$ of a positive linear Lebesgue measure such that the function
$Q$ is integrable over $\widetilde{S}(p_0, r)$ for every $r\in E_1.$
Then the family $\frak{Sob}^{*1\delta}_{K, Q}(D, D_*)$ consisting of
all homeomorphisms $f:D\rightarrow D_*,$ $f=g^{\,-1},$ $g\in
\frak{Sob^{-1}}^{\delta}_{K, Q}(D, D_*),$ is uniformly light on $K;$
moreover, $\frak{Sob}^{*1\delta}_{K, Q}(D, D_*)$ is equi-uniform on
$K.$

\medskip
Moreover, if (\ref{eq3*!!A})--(\ref{eq8BC}) hold for any $p_0\in D,$
any $0<r_1<r_2<\sup\limits_{p\in D}\widetilde{h}(p, p_0),$ any $f\in
\frak{Sob}^{*1\delta}_{K, Q}(D, D_*)$ and $Q\in L^1(D),$ then there
exist $C_1>0$ such that
$$
\widetilde{h}_*(f(x), f(y))\geqslant
C_1\cdot\exp\left\{-\frac{2\Vert Q\Vert_1}{C{(\widetilde{h}(x,
y))}^2}\right\}
$$
for all $x, y\in K,$ $x\ne y,$ and every $f\in
\frak{Sob}^{*1\delta}_{K, Q}(D, D_*),$ where $C$ is a constant from
(\ref{eq3C}) for $n=2$ and $X={\Bbb C}.$ }
\end{theorem}

\begin{corollary}\label{cor1_2}
{\, The statement of Theorem~\ref{th1_2} remains true if we replace
the assumption on $Q$ mentioned above concerning the integrability
of $Q$ over circles $\widetilde{S}(p_0, r)$ by a simpler condition:
$Q\in L_{\rm loc}^1(D).$}
\end{corollary}

\begin{corollary}\label{cor2_2}{\bf \,(An analogue of Koebe-Bloch
theorem)}.~{ Under assumptions of Theorem~\ref{th1_2}, the family
$\frak{Sob}^{*1\delta}_{K, Q}(D, D_*)$ is uniformly open on $K.$ In
particular, the latter is true if we replace the assumption on $Q$
mentioned in Theorem~\ref{th1_2} concerning the integrability of $Q$
over circles $\widetilde{S}(p_0, r)$ by a simpler condition: $Q\in
L_{\rm loc}^1(D).$}
\end{corollary}

\begin{theorem}\label{th2_2}
{\,Let $D$ and $D_*$ be domains in ${\Bbb S}$ and ${\Bbb S}_*,$
respectively, while $D$ has at least two different boundary points
$a$ and $b$ and, besides that, $\overline{D_*}$ is a compactum in
${\Bbb S}_*.$ Assume that, for each point $p_0\in D$ there is
$r_0=r_0(p_0)>0,$ $r_0<\sup\limits_{p\in D}\widetilde{h}(p, p_0),$
such that for any $0<r_1<r_2<r_0$ there is a set $E_1\subset[r_1,
r_2]$ of a positive linear Lebesgue measure such that the function
$Q$ is integrable over $\widetilde{S}(p_0, r)$ for every $r\in E_1.$
Then the family $\frak{Sob}^{*2\delta}_{A, Q}(D, D_*),$ consisting
of all homeomorphisms $f:D\rightarrow D_*,$ $f=g^{\,-1},$ $g\in
\frak{Sob}^{-2\delta}_{K, Q}(D, D_*),$ is uniformly light on every
compactum $K\subset D$ moreover, $\frak{Sob}^{*2\delta}_{A, Q}(D,
D_*)$ is equi-uniform on $K.$

\medskip
Moreover, if (\ref{eq3*!!A})--(\ref{eq8BC}) hold for any $p_0\in D,$
any $0<r_1<r_2<\sup\limits_{p\in D}\widetilde{h}(p, p_0),$ any $f\in
\frak{Sob}^{*2\delta}_{A, Q}(D, D_*)$ and $Q\in L^1(D),$ then there
exist $C_1>0$ such that the relation~(\ref{eq1D}) holds
for all $x, y\in K,$ $x\ne y,$ and every $f\in
\frak{Sob}^{*2\delta}_{A, Q}(D, D_*),$ where $C$ is a constant from
(\ref{eq3C}) for $n=2$ and $X={\Bbb C}.$ }
\end{theorem}

\begin{corollary}\label{cor6_2}
{\, The statement of Theorem~\ref{th2_2} remains true if we replace
the assumption on $Q$ mentioned above concerning the integrability
of $Q$ over circles $\widetilde{S}(p_0, r)$ by a simpler condition:
$Q\in L_{\rm loc}^1(D).$}
\end{corollary}

\begin{corollary}\label{cor7_2}{\bf \,(An analogue of Koebe-Bloch
theorem)}.~{ Under assumptions of Theorem~\ref{th2_2}, the family
$\frak{Sob}^{*2\delta}_{A, Q}(D, D_*)$ is uniformly open on every
compactum $K\subset D.$ In particular, the latter is true if we
replace the assumption on $Q$ mentioned in Theorem~\ref{th2_2}
concerning the integrability of $Q$ over circles $\widetilde{S}(p_0,
r)$ by a simpler condition: $Q\in L_{\rm loc}^1(D).$}
\end{corollary}

\medskip\medskip\medskip
{\bf Declarations.}

\medskip
{\bf Conflict of interest.} The authors have no financial or
proprietary interests in any material discussed in this article.

\medskip
{\bf Funding.} The work was supported by the National Research
Foundation of Ukraine (Project ``Analogues of Carath\'{e}odory and
Koebe-Bloch theorems for Orlycz-Sobolev classes'', Project number
2025.02/0010).

\medskip
{\bf Data availability statement.} All necessary data are included
in the paper.

\medskip
{\bf Contribution.} The authors have equal contribution to writing
the article.


\medskip
{\bf \noindent Miodrag Mateljevi\'{c}} \\ University of Belgrade,
Faculty
of Mathematics \\
16 Studentski trg, P.O. Box 550\\
11001 Belgrade, SERBIA\\
miodrag@matf.bg.ac.rs

\medskip
{\bf \noindent Evgeny Sevost'yanov} \\
{\bf 1.} Zhytomyr Ivan Franko State University,  \\
40 Velyka Berdychivska Str., 10 008  Zhytomyr, UKRAINE \\
{\bf 2.} Institute of Applied Mathematics and Mechanics\\
of NAS of Ukraine, \\
19 Henerala Batyuka Str., 84 116 Slov'yansk,  UKRAINE\\
esevostyanov2009@gmail.com

\end{document}